\newif\ifimfnum
\imfnumfalse

\documentclass[
\ifimfnum
english,twoside,abstracton,11pt,
\else
twoside,
\fi
]{\ifimfnum%
imfnumpreprint%
\else%
amsart%
\fi}

\ifimfnum
\preprintyear{2011}
\preprintnumber{1}
\else
\usepackage{graphicx}
\usepackage[pagebackref]{hyperref}
\fi

\newif\iffinal
\finaltrue

\usepackage{pdfsync}

\iffinal
\else
\usepackage[notref,notcite]{showkeys,rotating}

\fi

\usepackage{unixode}

\usepackage{amsmath,amsthm,amssymb}

\newcommand*{\mat}[1]{\mathsf{#1}}
\renewcommand*{\vec}[1]{\mathsf{#1}}
\newcommand{\NN}{\mathbf{N}}
\newcommand*{\bracket}[2]{\left\langle #1, #2\right\rangle}
\newcommand{\dd}{\mathrm d}
\newcommand*{\inv}{^{-1}}

\DeclareMathOperator{\im}{Im}

\newcommand{\Mv}{M}
\newcommand{\Vv}{V}
\newcommand{\Wv}{W}
\newcommand{\Nv}{N}
\newcommand{\Zv}{Z}
\newcommand{\Kv}{K}
\newcommand{\Cv}{C}
\newcommand{\Dv}{D}

\newcommand*{\red}[1]{^{(#1)}}
\newcommand*{\reds}{^{(∞)*}}

\DeclareMathOperator{\Span}{\mathrm{span}}
\newcommand{\coker}{\operatorname{coker}}

\DeclareMathOperator{\ind}{ind}

\newcommand{\xv}{\vec{x}}

\newcommand{\E}{\mat{E}}
\newcommand{\A}{\mat{A}}
\newcommand{\Pm}{\mat{P}}
\newcommand{\Qm}{\mat{Q}}
\newcommand{\Id}{\mathbb{I}}
\newcommand\transpose{{\!\scriptstyle\mathsf T}} 	

\newenvironment{eqsys}{\left\{\begin{aligned}}{\end{aligned}\right.}

\newcommand*{\cdom}[1]{_{(\E,\A)#1}}
\newcommand{\EA}{\cdom{}}

\newcommand{\Bb}{\mathcal{B}}
\newcommand{\Sb}{\mathcal{S}}

\newcommand{\DM}{∆\Mv}
\newcommand{\DV}{∆\Vv}

\newcommand{\CC}{\mathbf{C}}

\newcommand{\redss}{^{(∞)*(∞)}}
\newcommand{\redsss}{^{(∞)*(∞)*}}

\newcommand{\Rm}{\mat{R}}
\newcommand{\Eb}{\overline{\E}}
\newcommand{\Ab}{\overline{\A}}

\usepackage{xcolor}

\usepackage{aliascnt}
\iffinal
\colorlet{linkcolour}{black}
\colorlet{urlcolour}{black}
\else
\colorlet{linkcolour}{blue!50!black}
\colorlet{urlcolour}{magenta}
\fi
\hypersetup{
pdfauthor={Olivier Verdier},
colorlinks=true,
linkcolor=linkcolour, citecolor=linkcolour,
urlcolor=urlcolour,
}
\usepackage{enumerate}

\definecolor{weakorange}{rgb}{1,.9,.2}
\usepackage[textsize=footnotesize,backgroundcolor=weakorange,colorinlistoftodos,
]{todonotes}
\newcommand{\ltodo}[1]{\todo[backgroundcolor=green!40]{#1}}

\iffinal
\renewcommand{\todo}[1]{}
\fi

\usetikzlibrary{matrix}
\tikzstyle{commdiag}=[matrix of math nodes, row sep=3em, column sep=2.5em, text height=1.5ex, text depth=0.25ex]
\tikzstyle{exseq}=[commdiag, column sep=2em]

\tikzset{>=stealth}

\newcommand{\newtheoremalias}[2]{%
\newaliascnt{#1}{theorem}%
\newtheorem{#1}[#1]{#2}%
\aliascntresetthe{#1}%
}

\theoremstyle{plain}
\newtheorem{theorem}{Theorem}[section]

\newtheoremalias{proposition}{Proposition}

\newtheoremalias{corollary}{Corollary}

\newtheoremalias{lemma}{Lemma}

\theoremstyle{definition}
\newtheoremalias{definition}{Definition}

\newtheoremalias{notation}{Notation}

\theoremstyle{remark}
\newtheoremalias{example}{Example}

\newtheoremalias{remarkth}{Remark}

\newenvironment{remark}{\begin{remarkth}}{\end{remarkth}}

\newenvironment{pr}[1][\proofname]{
\begin{proof}[#1]}{\end{proof}
}

\newenvironment{thmenumerate}{\begin{enumerate}[{\normalfont (i)}]}{\end{enumerate}}

\newcommand{\picturesfolder}{.}
\newcommand{\matrixfigure}[4][.6]{\begin{figure}
\begin{tikzpicture}
	\coordinate (mat) at (0,0);
	\coordinate (leg) at (9,-1);
	\node[anchor=north west] at (mat) {\includegraphics[width=#1\textwidth]{\picturesfolder/#3}};
	\node[anchor=north west] at (leg) {\includegraphics[scale=.7]{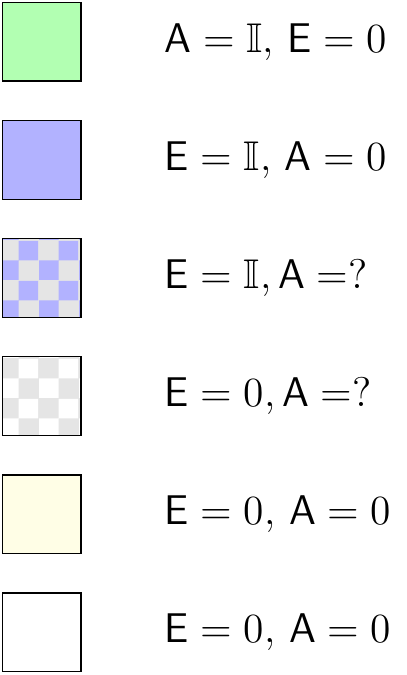}};
\end{tikzpicture}
\caption{#4}
\label{#2}
\end{figure}}

\newcommand*{\alert}[1]{\textit{\textbf{#1}}}

\title{Reduction and Normal Forms of Matrix Pencils}
\author{Olivier Verdier}

\ifimfnum

\fi

\begin{document}
\ifimfnum
	\pagestyle{empty}
	\cleardoublepage
	\pagestyle{plain}

\iffinal
\makeimfnumtitle
\fi

\else

\address{Department of Mathematical Sciences, NTNU, 7491 Trondheim, Norway}
\maketitle
\fi
\begin{abstract}
Matrix pencils, or pairs of matrices, may be used in a variety of applications.
In particular, a pair of matrices $(\E,\A)$ may be interpreted as the differential equation $\E x' + \A x = 0$.
Such an equation is invariant by changes of variables, or linear combination of the equations.
This change of variables or equations is associated to a group action.
The invariants corresponding to this group action are well known, namely the Kronecker indices and divisors.
Similarly, for another group action corresponding to the weak equivalence, a complete set of invariants is also known, among others the strangeness.

We show how to define those invariants in a directly invariant fashion, i.e. without using a basis or an extra Euclidean structure.
To this end, we will define a reduction process which produces a new system out of the original one.
The various invariants may then be defined from operators related to the repeated application of the reduction process.
We then show the relation between the invariants and the reduced subspace dimensions, and the relation with the regular pencil condition.
This is all done using invariant tools only.

Making special choices of basis then allows to construct the Kronecker canonical form.
In a related manner, we construct the strangeness canonical form associated to weak equivalence.
\end{abstract}


\iffinal
\else
\listoftodos
\fi


\section{Introduction}

\todo{Boast about relation defects/regpencil/invertibility}

\subsection{Equivalence}

The primary study of this paper is that of \emph{pairs of matrices}, also called \emph{matrix pencils}.
In other words, we study pairs of operators $(\E,\A)$ both acting from a finite dimensional vector spaces $\Mv$ to a finite dimensional vector space $\Vv$.

A typical example we have in mind is the linear differential equation
\begin{equation}
\label{eq_proto_diffeq}
\E \frac{\dd x}{\dd t} + \A x = 0
.
\end{equation}
Such a model is clearly invariant by changes of variable, or by changing the order of the equations.
More precisely, it is invariant by simultaneous \emph{equivalence transformation} of the operators $\E$ and $\A$.
The corresponding equivalence relation is the following: two pairs of operators $(\E_1,\A_1)$ and $(\E_2,\A_2)$ will be considered equivalent if there exists invertible operators $\mat{P}$ and $\mat{Q}$, operating on $\Mv$ and $\Vv$ respectively, such that
\begin{equation}
\label{eq_equiv_rel}
\begin{split}
\E_2 = \mat{P}\E_1\mat{Q}
,\\ 
\A_2 = \mat{P}\A_1\mat{Q} 
.
\end{split}
\end{equation}
This equivalence relation is associated to a group, which is simply $GL(\Mv)\times GL(\Vv)$.
This is called \emph{strong equivalence} in \cite{Gantmacher}.
We are interested in properties which are invariant with respect to that group action on the matrix pencil.
In other words, we are interested in quantities that label the orbit of the group action.

In fact, a complete set of invariants and a canonical form have been known since the works of \cite{Weierstrass} and \cite{Kronecker}.
Modern versions of those proofs may be found in \cite[§~XII.4]{Gantmacher} and in \cite[§~A.7]{Gohberg}.
The primary tool for obtaining those invariants is the Jordan canonical form.
For that reason, those proofs are impossible to extend to nearby cases, for example to the infinite dimensional case, or to the parameter dependent case, not to mention the numerical difficulties associated with the computation of the Jordan canonical form.

As a result, alternative proof techniques were developed, most notably in \cite{vanDooren} and \cite{Wilkinson}.
Those authors observed indeed that using a Jordan canonical form is not suitable to compute the invariants other than the Jordan invariants, i.e., the Kronecker indices and the ``infinite elementary divisors'' \cite{vanDooren}.
The idea is to transform the pair of matrices into a form which exhibits all the invariants but is not a canonical form.
Those forms are known under the names of generalized Schur-staircase form, or GUPTRI (Generalized Upper Triangular Form).
We refer to \cite[§4.1]{Johansson} and \cite{matrixpencils} for more references on those algorithms.

Our approach is similar, although with a shift of focus towards the underlying algebraic structures as opposed to the algorithmic aspects.
In particular, we attempt to define the invariants from the dimensions of subspaces which are themselves invariants with respect to the equivalence relation at hand.
The advantage of our approach is that a great deal of results are automatically independent of the choice of a basis, or any other structure (like a Euclidean structure).

\subsection{Invariants}

A matrix pencil, when considered as a differential equation \eqref{eq_proto_diffeq}, may be decomposed in an intrinsic ordinary differential equation, and an extra structure.
We will call the invariants of the underlying ordinary differential equation the \emph{dynamical} invariants, and we will call the remaining invariants the \emph{non-dynamical invariants}.
In the parlance of the Kronecker decomposition theorem as presented in \cite{Gantmacher}, the dynamical invariants would be the finite elementary divisors (essentially a Jordan form), whereas the non-dynamical invariants would be the infinite elementary divisors along with the row and column minimal indices.

The dynamical invariants, i.e., the invariants of the intrinsic differential equations boil down to the Jordan invariant associated to similarity transformations, and are therefore of less interest to us.
We will thus mostly focus on the non-dynamical invariants, which appear only when $\E$ is not invertible.
Those invariants are well-known in control theory, and in the study of differential algebraic equations.
In control theory, such invariants are the controllability and observability indices (\cite[§~6.3]{Kailath}), for differential algebraic equations (DAE),
in the case of a \emph{regular pencil} (see \autoref{sec_regular_pencil}), the most used non-dynamical invariant is the \emph{index} (\cite[VII.1]{Hairer2}).

Our goal is to define the non-dynamical invariants in a invariant manner, without any other structure than the linear algebraic structure.

\subsection{Reduction}

The crucial tool to the study of the invariants of a pencil is the concept of \emph{reduction}, which we define precisely in \autoref{sec_reduction}.

This concept was gradually developed, under various names, or no name at all, first in \cite{Wong} for the study of regular pencils, then in \cite[§4]{Wilkinson} and \cite{vanDooren} to prove the Kronecker decomposition theorem.
It is also related to the \emph{geometric reduction} of nonlinear implicit differential equations as described in \cite{Reich} or \cite{Rabier-Rheinboldt}. 
In the linear case, those coincide with the observation reduction, as shown in \cite{thesis}.
It is also equivalent to the algorithm of prolongation of ordinary differential equation in the formal theory of differential equations, as shown in \cite{Reid}.

The reduction procedure is an operation that, out of a pair of operators $(\E,\A)$, creates a new, smaller one $(\E',\A')$.
``Smaller'' is in the sense that the reduced operators $\E'$ and $\A'$ are restrictions of $\E$ and $\A$ on subspaces of $\Mv$ and $\Vv$, defined by $\Vv':=\E\Mv$ and $\Mv':=\A\inv\Vv'$.

This process of reduction is iterated, producing systems $(\E\red{k},\A\red{k})$ and subspaces $\Mv\red{k}$ and $\Vv\red{k}$.
This process ultimately stops, and we will call the number of steps before it stops the \emph{index}.
\todo{bias in that name: obs/cont index}
When the process stops, the system which is produced, denoted by $(\E\red{∞},\A\red{∞})$, is such that $\E\red{∞}$ is surjective.
After running the reduction algorithm once more on the dual of that reduced system, i.e., on $(\E\reds,\A\reds)$, one obtains an isolated system $(\E\redss,\A\redss)$ such that $\E\redss$ is now invertible.

At each step of the reduction, some information from the original system is lost.
That information is encoded by integers called ``defects''.
Those defects are of three kinds: $α$, $β^+$ and $β^-$.
The defect $α_1$ is defined as the dimension of the kernel of  $\E$, regarded as a quotient operator from $\Mv/\Mv'$ to $\Vv'/\Vv''$. 
The defect $β^+_1$ is defined as the dimension of the cokernel of $\A$, regarded as a quotient operator from $\Mv/\Mv'$ to $\Vv/\Vv'$.
The iterated reduction then generates the sequences of defects $α_k$ and $β^+_k$.
The defects $β^-_k$ are defined as the $β^+$ defects of the system $(\E\reds,\A\reds)$.

Using those subspaces, defined in an invariant manner, we are able to show the following facts:
\begin{itemize}
	\item the operator $\E$ is invertible if and only if all the defects vanish
	\item the pair $(\E,\A)$ is a regular pencil if and only if the $β^+$ and $β^-$ defects vanish
	\item we show that the invariants defined in \cite{Mehrmann}, like the strangeness, may also be defined directly in an invariant manner, i.e., without using any extra structure or basis
\end{itemize}

We also show that the defects and the system $(\E\redss,\A\redss)$ completely characterize the equivalence class corresponding to the equivalence relation \eqref{eq_equiv_rel}.
\begin{itemize}
	\item the defects are related to the Kronecker indices
	\item the invariants defined in \cite{Mehrmann} may be used to construct a corresponding canonical form for weak equivalence: this connects the approaches of \cite{Wilkinson} and \cite{Mehrmann}
	\item using the relation with the Kronecker decomposition theorem, we show that the defects of the dual system $(\E^*,\A^*)$ are related to those of $(\E,\A)$ by switching the $β^+$ and $β^-$ defects.
\end{itemize}

\subsection{Outline}

The layout of the paper is as follows.

In the first part, \autoref{sec_reduction} and \autoref{sec_defects}, we show how to derive the non-dynamical invariants.
In \autoref{sec_reduction} we define the reduction procedure. 
In \autoref{sec_defects} we define the defects of a system, and study their properties. 
In particular, we give an original proof of the relation between the property of a pencil to be regular, and the presence of some of the defects.

In the second part, we show that the invariants obtain in the first part, namely the defects, supplemented by a Jordan structure, are the only invariants of the pair of matrices with respect to equivalence.
Most of the results in this part are already in \cite{Wilkinson} and \cite{vanDooren}.
In \autoref{sec_coupling} we prove the basic lemmas needed to construct canonical forms.
In \autoref{sec_strangeness}, we show how to use those tools to construct a canonical form with respect to weak equivalence
In \autoref{sec_decomposition} and \autoref{sec_decomposition_bis} we show that the defects determine a complete canonical form.
In \autoref{sec_kronecker} we study the relation with the existing Kronecker canonical form.

\section{System Reduction}
\label{sec_reduction}

\subsection{Setting}

\begin{definition}
\label{def_system}
We will call a pair of linear operators $(\E,\A)$ a \alert{linear system}, or simply a \alert{system}, if $\E$ and $\A$ have the same domain and codomain, both of finite dimension.
\end{definition}

Given a system $(\E,\A)$, we will denote the common \alert{domain} of $\E$ and $\A$ by 
$\Mv\EA$ and the common \emph{codomain} of $\E$ and $\A$ by $\Vv\EA$, so 
a system $(\E,\A)$ may be represented as

\[
\E,\A:\Mv\EA\longrightarrow \Vv\EA
.
\]


\subsection{Reduced spaces}
\label{red_space_sec}

The idea behind the reduction of a linear system $(\E,\A)$ is to ``disentangle'' the spaces associated with the operators $\E$ and $\A$. 
The strategy pursued is to try and make the operator $\E$ surjective, by successive reduction steps.
In order to achieve this, we have to describe the lack of surjectivity  of $\E$, first independently of $\A$, which leads to the definition of the subspace
\[\Vv' := \E\Mv
.
\]
The next step is now to describe the lack of surjectivity of $\E$, \emph{with respect to $\A$}, which we measure using the subspace
\[\Mv' := \A\inv\Vv'
.
\]

\begin{remark}
Those definitions make sense when considering the differential equation
\[\E x' + \A x = 0
.
\]
Notice that any suitable initial condition for this equation must be in $\Mv'$. If the initial condition is not in $\Mv'$, there cannot be any solution stemming from that initial condition.
\end{remark}


Let us put those definitions together:

\begin{definition}
\label{space_red_def}
Given a linear system $(\E,\A)$ we define its \alert{reduced codomain} {$\Vv'\EA$}{Reduced codomain} as
\[\Vv'\EA := \E\Mv\EA
,
\]
and its \alert{reduced domain} $\Mv'\EA$ as
\[\Mv'\EA := \A\inv\Vv'\EA = \bigl\{\xv\in\Mv\EA:\ \A \xv \in \Vv'\EA \bigr\} 
.
\]
\end{definition}

\begin{remark}
	We will often drop the dependency on the system $(\E,\A)$, and simply write $\Mv$, $\Mv'$, $\Vv$ and $\Vv'$ when the context is clear enough.
\end{remark}


\begin{remark}
As explained in \cite[§5.1]{thesis}, the reduction of \autoref{space_red_def} corresponds to the non-linear reduction of general systems of differential equations with constraints. 
The study of differential equations is also the point of departure in \cite{Wilkinson}.
\end{remark}

\begin{remark}
One of the first occurrence of the definition of that subspace $\Mv'$ seems to be in \cite[Lemma 2.1]{Wong}. It is used to study systems which are regular pencils (see \autoref{def_reg_pencil}).

Another explicit definition is to be found in \cite[§7]{Reich-report}, although with a different purpose than ours, namely the study of linear, time-varying differential algebraic equations of index one.

\end{remark}


\newcommand{\restr}{|_{\Mv'}}

\subsection{System Reduction}
\label{system_reduction_sec}

The subspaces $\Mv\red{k}$ and $\Vv\red{k}$ defined in \autoref{space_red_def} allow for defining a new system.
This procedure will be called ``reduction''.

\begin{proposition}
\label{prop_red_ops}
	Given a system $(\E,\A)$, the operators $\E'$ and $\A'$ are uniquely defined by the following commuting diagram.
\begin{center}
\begin{tikzpicture}
\matrix(m) [commdiag, column sep=3.5em]
{\Mv & \Vv\\
\Mv' & \Vv' \\};
\path[->]
(m-1-1) edge node[auto] {$\E,\A$} (m-1-2)
(m-2-1) edge node[auto] {$\E',\A'$} (m-2-2)
(m-2-1) edge (m-1-1)
(m-2-2) edge (m-1-2);
\end{tikzpicture}
\end{center}
The vertical arrows are canonical injection from a subspace into the ambient space.

The operators $\E'$ and $\A'$ build up a new system $(\E,\A)'$ which we call the \alert{reduced system}, and is defined by
\[(\E,\A)' := (\E', \A')
.
\]
\end{proposition}
\begin{pr}
The proof rests on the observation that
\[\E\Mv'\EA \subset \Vv'\EA\qquad\text{and}\qquad \A\Mv'\EA \subset \Vv'\EA
.
\]
\end{pr}

\begin{remark}
Consider the category which objects are vector spaces and arrows are systems as defined in \autoref{def_system}.
The reduction operation, denoted by a prime, is an \emph{endofunctor} in this category, i.e., a functor from that category to itself.
\end{remark}

As we mentioned in the beginning of \autoref{red_space_sec}, our goal is to obtain a reduced system such that $\E$ is surjective.
It is only part of a general strategy to obtain a reduced system where $\E$ is invertible.
It is therefore important that the reduction algorithm does not alter the injectivity of $\E$.
We observe that this is indeed the case.

\begin{proposition}
\label{prop_Ered_inj}
If, in a system $(\E,\A)$, $\E$ is injective, then $\E'$ is also injective.
\end{proposition}
\begin{pr}
It is a consequence of the observation that
\[\ker\E' \subset \ker\E
.
\]
\end{pr}

The pendant of that observation is the equally simple observation regarding the kernel of the operator $\A$ with respect to the reduced space $\Mv'$:

\begin{proposition}
\label{kerA_sub_Mprim_prop}
Given a system $(\E,\A)$, the null-space of $\A$ is included in $\Mv'\EA$, i.e.,
\[\ker\A \subset \Mv'\EA
.
\]
\end{proposition}

\subsection{Iterated Reduction}



We may iterate the reduction process described in \autoref{system_reduction_sec} on the new system $(\E,\A)'$. 
This leads to a sequence of systems $\{(\E,\A)\red{k}\}_{k\in\NN}$
which is defined recursively as follows.
\begin{definition}
\label{recsys_red_def}
The iterated of the reduction of a system $(\E,\A)$ are defined recursively by
\[(\E\red{k+1},\A\red{k+1}) := (\E\red{k},\A\red{k})'
,
\qquad \forall k≥0
,
\]
and
\[(\E\red{0},\A\red{0}) := (\E,\A)
.
\]
\end{definition}

We will make use of the straightforward notation, for $k\in\NN$.
\begin{equation}
\label{space_red_notation}
\begin{aligned}
\Mv\red{k}\EA := \Mv\cdom{\red{k}},\\
\Vv\red{k}\EA := \Vv\cdom{\red{k}}.
\end{aligned}
\end{equation}


The reduced operators $\E\red{k}$ and $\A\red{k}$ are essentially restrictions of the original operators $\E$ and $\A$, so we may rewrite the definition of the iterated reduced subspaces $\Mv\red{k}$ and $\Vv\red{k}$.
\begin{proposition}
\label{red_alt_def_prop}
For a system $(\E,\A)$ the following assertions hold for any integer $k≥0$:
\[\forall x\in\Mv\red{k}\quad \E\red{k}x = \E x \quad \A\red{k}x = \A x
,
\]
\[\Vv\red{k+1} = \E \Mv\red{k}
,
\]
\[\Mv\red{k+1} = \bigl\{x\in\Mv\red{k}:\ \A x \in \Vv\red{k+1}\bigr\}
.
\]
\end{proposition}
\begin{pr}
The proof is a simple verification by induction on $k$.
\end{pr}

\subsection{Totally Reduced Systems}

As we shall notice in \autoref{sec_index}, the repeated operation of reduction transforms a system into one which cannot be reduced anymore, or rather, for which the reduction does not create a new system.
We call such systems ``totally reduced'':

\begin{definition}
We will say that a system $(\E,\A)$ is \alert{totally reduced} if
\[(\E,\A)' = (\E,\A)
.
\]
\end{definition}

A practical characterisation of a totally reduced system is that $\Vv' = \Vv$. The verification is straightforward.

\begin{proposition}
\label{prop_tot_red_V}
A system $(\E,\A)$ is totally reduced if and only if
\[\Vv'\EA = \Vv\EA
.
\]
\end{proposition}

\subsection{Almost Reduced System}

\begin{definition}
\label{def_almost_red}
We will say that a system $(\E,\A)$ is \alert{almost reduced} if
\[\Mv'\EA = \Mv\EA
.
\]
\end{definition}

The chosen vocabulary is supported by the following facts:
\begin{itemize}
	\item a totally reduced system is also almost reduced, which follows from \autoref{space_red_def}.
	\item a system which is almost reduced will be totally reduced at the next step of the reduction, since by \autoref{red_alt_def_prop}: $\Vv'' = \E\Mv' = \E\Mv = \Vv'$.
\end{itemize}

\begin{remark}
In the situation of a reduced system which is almost but not totally reduced,
 the following subspace sequences
\begin{equation*}
\label{eq_nested_spaces}
\begin{aligned}\Mv\red{n+1}=&\Mv\red{n}\subset \ldots \subset \Mv'' \subset \Mv' \subset \Mv \\
\Vv\red{n+2}=\Vv\red{n+1}\subset &\Vv\red{n}\subset \ldots \subset \Vv'' \subset \Vv' \subset \Vv
\end{aligned}
\end{equation*}
would be produced.

A concrete example where this happens is when $\A=0$ and $\E$ is not surjective.
It is clear that $\Mv' = \Mv$ but $\Vv' \subsetneq \Vv$.
The corresponding system is thus almost reduced but not totally reduced.

\end{remark}

\subsection{Index}
\label{sec_index}

The reduction procedure produces decreasing sequences of subspaces.
When both sequences stall, the system is totally reduced.
The number of reduction steps needed to transform a system into a totally reduced one is called the \emph{index} of the system $(\E,\A)$:

\begin{definition}
\label{def_index}
The smallest integer $n\in\NN$  for which the system $(\E\red{n},\A\red{n})$ is totally reduced
is called the \alert{index} of the system $(\E,\A)$. 

We will use the following notation for the index of the system $(\E,\A)$:
\[\ind\EA := \min \bigl\{ n\in\NN:\ (\E,\A)\red{n+1} = (\E,\A)\red{n} \bigr\}
.
\]
\end{definition}

\begin{remark}
The index is \emph{always} a finite integer\footnote{as opposed to the differentiation index, which is infinite in the non-regular pencil case; see, e.g., \cite[§ VII.1]{Hairer2}.}, and the reduced system $(\E',\A')$ has an index dropped by one, i.e,
\[\ind\cdom{'} = \ind\EA - 1
.
\]

Those observations will be used repeatedly to prove statements by induction on the index (e.g., in \autoref{prop_reg_pencil}, \autoref{thm_linear_decomposition} and \autoref{thm_NL_basis}).
\end{remark}

\begin{remark}
\label{rem_def_index}
Using \autoref{prop_tot_red_V} we observe that
\[\ind\EA = \min \bigl\{ n\in\NN:\ \Vv\red{n+1} = \Vv\red{n} \bigr\}
.
\]
\end{remark}

\begin{remark}
The index defined in \autoref{def_index} is closely related to the geometric index defined in \cite{Reich}, \cite{Rabier-Rheinboldt} or \cite[§5.1]{thesis}.
In fact, the geometric index would be the first integer $n$ such that the system $(\E,\A)\red{n}$ is \emph{almost} reduced (\autoref{def_almost_red}).
As we shall see in \autoref{prop_ind_no_beta} and \autoref{prop_reg_pencil}, this minor difference is only relevant for singular pencils.
\end{remark}

\subsection{Totally Reduced System}

\begin{definition}
\label{def_totred_sys}
For a system $(\E,\A)$ of index $n=\ind\EA$ we define the \alert{totally reduced system} as
\[(\E\red{∞},\A\red{∞}) := (\E\red{n},\A\red{n})
.
\]
\end{definition}

\begin{remark}
We could simply have defined, say $\E\red{∞}$ by the limit of the  sequence of operators $\E\red{k}$ (because this sequence eventually stalls), which explains the notation ``$∞$''.
\end{remark}

We pointed out in \autoref{red_space_sec} that the idea behind the reduction procedure was to lead to a system where $\E$ is surjective.
The reduction algorithm indeed achieves this goal:

\begin{proposition}
\label{prop_Eredinf_surj}
The totally reduced operator $\E\red{∞}$ is surjective.
\end{proposition}
\begin{pr}
The system $(\E\red{∞},\A\red{∞})$ is totally reduced so we may use \autoref{prop_tot_red_V} to conclude that $\E\red{∞} \Mv\red{∞} = (\Vv\red{∞})' = \Vv\red{∞}$, so $\E\red{∞}$ is surjective.
\end{pr}

\section{Defects}
\label{sec_defects}

\subsection{Quotient Operators}

At each step of the reduction some information is lost, by passing from the original system to the reduced one.
We capture that information loss by two quotient operators defined on the quotient space $\Mv/\Mv'$.

\begin{proposition}
\label{prop_quotient_op}
The following commutating diagrams uniquely define the quotient operators $[\A]$ and $[\E]$  (the vertical arrows are the natural projections on a quotient space).
\begin{center}
\begin{tikzpicture}
\matrix(m) [commdiag]
{\Mv & \Vv\\
\Mv/\Mv' & \Vv/\Vv' \\};
\path[->]
(m-1-1) edge node[auto] {$\A$} (m-1-2)
(m-2-1) edge node[auto] {$[\A]$} (m-2-2)
(m-1-1) edge (m-2-1)
(m-1-2) edge (m-2-2);

\begin{scope}[xshift=30ex]
\matrix(m) [commdiag]
{\Mv & \Vv'\\
\Mv/\Mv' & \Vv'/\Vv'' \\};
\path[->]
(m-1-1) edge node[auto] {$\E$} (m-1-2)
(m-2-1) edge node[auto] {$[\E]$} (m-2-2)
(m-1-1) edge (m-2-1)
(m-1-2) edge (m-2-2);
\end{scope}
\end{tikzpicture}
\end{center}

Moreover, $[\A]$ is injective and $[\E]$ is surjective. 
\end{proposition}

\begin{pr}
The quotient operators $[\A]$ and $[\E]$ are well defined because $\A\Mv' \subset \Vv'$ and $\E\Mv' \subset \Vv''$ (since in fact, $\E\Mv' = \Vv''$ by definition). 
$[\E]$ is surjective because $\E$ is surjective onto $\Vv'$ by definition of $\Vv'$. 
$[\A]$ is injective since, by definition of $\Mv'$,
\[\A x \in\Vv' \implies x\in\Mv'
.
\]
\end{pr}

\subsection{Constraint and Observation Defects}

Since $[\A]$ is injective and $[\E]$ is surjective, the information stemming from those operators are to be collected in the cokernel of $[\A]$ and the kernel of $[\E]$.
The dimension of those subspaces are important invariants of the system $(\E,\A)$ which we now precisely define.

\begin{definition}
\label{def_defects}

Let $[\E]$ and $[\A]$ be defined as in \autoref{prop_quotient_op}.
We measure the lack of surjectivity of $[\A]$ by the \alert{first observation defect} $β^+_1(\E,\A)$, defined as
\begin{equation}
\label{def_alpha_one}
	β^+_1(\E,\A) := \dim \coker [\A] 
\end{equation}
and the lack of injectivity of $[\E]$ by the \alert{first constraint defect} $α_1(\E,\A)$, defined as
\begin{equation}
\label{def_betap_one}
	α_1(\E,\A) := \dim \ker [\E]
.
\end{equation}
\end{definition}

Now we take advantage of the reduction procedure and define those defects recursively:
\begin{definition}
\label{def_red_defects}
The \alert{constraint defects} $α_{k}(\E,\A)$ of a system $(\E,\A)$ are defined for any integer $k≥1$ by
\[α_{k}(\E,\A) := α_{1}((\E,\A)\red{k-1})
.
\]
Similarly, the \alert{observation defects} $β_{k}^+(\E,\A)$ are defined for any integer $k≥1$ by
\[β^+_{k}(\E,\A) := β_{1}^+((\E,\A)\red{k-1})
.
\]
\end{definition}

\subsection{Control Defects}

There is another important kind of defect that will be needed.
In is obtained by considering the dual of the totally reduced system obtained after repeated reductions.
That totally reduced system $(\E\red{∞},\A\red{∞})$ is such that $\E\red{∞}$ is surjective, so $\E\reds$ is injective.
So what happens for a system $(\E,\A)$ such that $\E$ is injective?
It turns out that such a system has no constraint defects.

\begin{proposition}
\label{prop_Einj_noalpha}
Given a system $(\E,\A)$, if $\E$ is injective, then the system has no constraint defects, i.e., for all integer $k≥1$, $α_k(\E,\A) = 0$.
\end{proposition}
\begin{pr}
The proof proceeds by induction on the index.
\begin{enumerate}
\item
If the index is zero, then $\Mv'=\Mv$ and $\Vv'=\Vv$, so $\dim\ker[\E]=0$. 
\item
For a positive index, using \autoref{prop_Ered_inj} we may apply the induction hypothesis and deduce that $α_{k}(\E,\A)=0$ for $k≥2$.
\item
Now if $x+\Mv' \in\ker[\E]$ then $\E x\in\E\Mv'$. 
Since $\E$ is injective, this means that $x\in\Mv'$ and thus that $\ker[\E]=0$. We conclude that $α_1(\E,\A)=0$.
\end{enumerate}
\end{pr}

Let us introduce the notion of a dual system.
\begin{notation}
Given a system $(\E,\A)$ we define the \alert{dual system} $(\E,\A)^*$ by the pair of adjoint operators $(\E^*,\A^*)$, i.e.,
\[(\E,\A)^* := (\E^*,\A^*)
.
\]
\end{notation}

\begin{proposition}
\label{prop_alpha_reds_zero}
The dual $(\E\red{∞},\A\red{∞})^*$ of a totally reduced system has no constraint defects, i.e.,
\[α(\E\reds,\A\reds) = 0
.
\]
\end{proposition}
\begin{pr}
According to \autoref{prop_Eredinf_surj}, the operator $\E\red{∞}$ is surjective, so $\E\reds$ is injective, and we conclude using \autoref{prop_Einj_noalpha}.
\end{pr}

This suggests that another set of defects is given by the observation defects of the dual of the totally reduced system $(\E\red{∞},\A\red{∞})$.
\begin{definition}
\label{def_betam}
Given a system $(\E,\A)$, we define the \alert{control defects} $β^-_k(\E,\A)$ by
\begin{equation*}
\label{eq_def_betam}
β^-_k(\E,\A) := β^+_k(\E\reds,\A\reds)\qquad\forall k≥1.
\end{equation*}
\end{definition}

\subsection{Intrinsic Dynamical System}

The reduction procedure may thus be used once to obtain a totally reduced system, and may then be applied again to the dual of that totally reduced system.

Starting with a system $(\E,\A)$, we may completely reduce it to obtain the system $(\E\red{∞},\A\red{∞})$.
The operator $\E\reds$ is injective.
The adjoint system $(\E\reds,\A\reds)$ may be in turn completely reduced to obtain the system $(\E\redss,\A\redss)$.
Using \autoref{prop_Eredinf_surj} and \autoref{prop_Ered_inj}, we obtain the following result.

\begin{proposition}
\label{prop_Edyn_inv}
The operator $\E\redss$ is invertible.
\end{proposition}

Since the operator $\E\redss$ is invertible, its domain and co-domain have the same dimension. 
This dimension is the dimension of the intrinsic dynamics of the system.

\begin{definition}
\label{def_dyndim}
The \alert{dynamical dimension} $δ$ of the system $(\E,\A)$ is defined by the integer
\[δ := \dim\Mv\redss = \dim\Vv\redss
.
\]
\end{definition}

\begin{remark}
For a differential equation defined by the system $(\E,\A)$, the system 
\[(\E\redsss,\A\redsss)\]
 corresponds to the underlying differential equation.
In particular, the dynamical dimension $δ$ determines the degrees of freedom for the choice of the initial condition.
\end{remark}


\subsection{Dimensions of the Subspaces}

In order to study the relations existing between the defects and the various subspaces $\Mv\red{k}$ and $\Vv\red{k}$, we define the following spaces, which measure the difference of dimension between each successive reduction:

\begin{definition}
\label{def_DMDV}
Recalling \autoref{recsys_red_def}, for any integer $k≥1$ we define the spaces
\[\DM\red{k} := \Mv\red{k-1}/\Mv\red{k}
\]
and
\[\DV\red{k} := \Vv\red{k-1}/\Vv\red{k}
.
\]
\end{definition}

By definition of the defects in \autoref{def_defects} and using \autoref{prop_quotient_op}, one obtains the relations
\begin{equation}
\label{eq_DVDM}
\begin{aligned}
\dim \DM\red{k} &= \dim\DV\red{k+1} + α_k
,
\qquad \forall k≥1
,\\
\dim \DV\red{k} &= \dim \DM\red{k} + β_k^+
,
\qquad \forall k≥1
,
\end{aligned}
\end{equation}
 between the dimensions of the spaces defined in \autoref{def_DMDV} and the defects.

For any integer $k≥1$ this implies the inequalities
\[ \cdots \leq \dim\DM\red{k+1} ≤ \dim \DV\red{k+1} ≤ \dim \DM\red{k-1} ≤ \dim \DV\red{k}\leq \cdots
.
\]

\begin{remark}
	This is the same sequence of inequalities as in \cite[5.2]{Wilkinson}.
\end{remark}


In particular, the dimensions of the spaces $\DM\red{k}$ and $\DV\red{k}$ may be expressed using the constraint and observation defects.

\begin{lemma}
\label{lma_rel_defects_dim_DVDM}
For any integer $k≥1$, the dimensions of the spaces $\DM\red{k}$ and $\DV\red{k}$ are related to the defects by the identities
\[\begin{aligned}
\dim \DV\red{k} &= ∑_{j≥k}(α_j + β^+_j)%
,%
\\
\dim \DM\red{k} &= ∑_{j≥k}(α_{j} + β^+_{j+1})%
.%
\end{aligned}\]
\end{lemma}
\begin{pr}
Those identities follow from an induction based on \eqref{eq_DVDM} and the observation that the integers $\dim\DV\red{k}$ and $\dim\DM\red{k}$ are zero when $k$ is bigger than the index of the system.
\end{pr}

\begin{remark}
As we shall see in \autoref{thm_strangeness}, the quantity defined in \cite{Mehrmann} as the ``strangeness'' $s$ turns out to be the integer
\[
s = \dim \DV''
.
\]
Roughly speaking it expresses the number of constraints that, when differentiated, will help to reduce the system.

We may thus give the precise relation of the strangeness to the defects using \autoref{lma_rel_defects_dim_DVDM}, namely
\[s = \dim\DV'' = ∑_{k=2}^{∞} β_k^+ +  ∑_{k=2}^{∞} α_k
.
\]
\end{remark}

The dimensions of the spaces $\Mv$ and $\Vv$ may also be expressed from the defects and the dynamical dimension $δ$ (see \autoref{def_dyndim}).

\begin{proposition}
\label{prop_rel_defect_dims}
The dimensions of $\Mv$, $V$, the defects $α$, $β^+$ and $β^-$ and the dynamical dimension $δ$ are related by the formulae
\[\begin{aligned}
	\dim \Mv &= δ + ∑_{k≥1} k α_k + ∑_{k≥1} k β^-_k + ∑_{k≥1}kβ^+_{k+1}
,\\
	\dim \Vv &= δ + ∑_{k≥1} k α_k + ∑_{k≥1} k β^+_k + ∑_{k≥1}k β^-_{k+1}
.
\end{aligned}\]
\end{proposition}
\begin{pr}
First observe that since $\dim \Vv\red{k} = \dim\Vv\red{k+1}+\dim\DV\red{k+1}$ and $\dim\Mv\red{k} = \dim\Mv\red{k+1} + \dim\DM\red{k+1}$, we have
\[\dim\Mv = \dim\Mv\red{∞} + ∑_{k≥1} \dim\DM\red{k} \qquad \dim \Vv = \dim\Vv\red{∞} + ∑_{k≥1}\dim\DV\red{k}
.
\]
Using \autoref{lma_rel_defects_dim_DVDM} we obtain
\[\dim\Vv = \dim\Vv\red{∞} + ∑_{k≥1}k α_k + ∑_{k≥1}k β_k^+
,
\]
and
\[\dim\Mv = \dim\Mv\red{∞} + ∑_{k≥1}kα_k +∑_{k≥1}kβ_{k+1}^+
.
\]
Now using the observation of \autoref{prop_alpha_reds_zero} that $α(\E\reds,\A\reds) = 0$, along with \autoref{def_betam} of the defects $β^-$ and \autoref{def_dyndim} of the dynamical dimension $δ$ we readily obtain the result.
\end{pr}

\subsection{Relation with the Index}

The index is, as expected, a non-dynamical invariant.
More precisely, it is a function of the defects, as the following proposition shows:
\begin{proposition}
The {index} $\ind\EA$ (see \autoref{def_index}) of a linear system $(\E,\A)$ is given by
\[\ind\EA = \min \bigl\{n\in\NN:\ \forall k > n\quad α_k(\E,\A) =0\quad\text{and}\quad β^+_{k}(\E,\A) =0\bigr\}
.
\]
\end{proposition}
\begin{pr}
Following \autoref{rem_def_index}, the index fulfills
\[\ind\EA = \min_k \dim\DV\red{k+1}=0
.
\]
Using \autoref{lma_rel_defects_dim_DVDM} we thus obtain
\[\dim \DV\red{k} = 0 \iff α_j+β^+_{j} =0\quad\forall j ≥ k+1
,
\]
which proves the claim.
\end{pr}

In the case of a system without observation defects we obtain readily:
\begin{corollary}
\label{prop_ind_no_beta}
The {index} of a system $(\E,A)$ without observation defects (i.e., $β^+=0$) is the biggest index of non-zero constraint defects, i.e.,
\[\ind\EA = \min \bigl\{n\in\NN:\ \forall k>n\quad α_k(\E,\A) = 0\bigr\}
.
\]
\end{corollary}
\ltodo{add remark on index for DAEs?}

\subsection{Defects and Invertibility}

The choice of the name ``defect'' may seem overly negative, but those integers really measure how far this system is from a system where $\E$ is invertible.
This is the essence of the following proposition.

\begin{proposition}
For a given system $(\E,\A)$ the following statements are equivalent.
\begin{thmenumerate}
	\item All the defects $\alpha$, $\beta^+$ and $\beta^-$ are zero.
	\item The operator $\E$ is invertible.
\end{thmenumerate}
\end{proposition}

\begin{pr}
$\E$ is surjective if and only if $\DV' = 0$. By \autoref{lma_rel_defects_dim_DVDM}, that is equivalent to $α = β^+ = 0$.
Since $\E$ is invertible if and only if both $\E$ and $\E^*$ are surjective, we obtain the result using \autoref{def_betam}.
\end{pr}

\subsection{Regular Pencils}
\label{sec_regular_pencil}

A \emph{pencil} is a polynomial on a ring of matrices.
Since we are interested in pairs of matrices, our attention is restricted to first order polynomials, and to the property of such a polynomial to be \emph{regular}.

\begin{definition}
\label{def_reg_pencil}
The system $(\E,\A)$ is a \alert{regular pencil} if there exists $λ\in\CC$ such that $λ\E + \A$ is invertible.
\end{definition}

There is a remarkable relation between the property of being regular and the defects:

\begin{proposition}
\label{prop_reg_pencil}
The system $(\E,\A)$ is a regular pencil if and only if all the defects $β^+$ and $β^-$ are zero.
\end{proposition}

We need first a lemma to understand how the pencil regularity property may be lost during the reduction.

\begin{lemma}
\label{lma_regular_pencil}
The system $(\E,\A)$ is a regular pencil if and only if both the following properties hold:
\begin{thmenumerate}
	\item $β_1^+(\E,\A)=0$
	\item The reduced system $(\E,\A)'$ is a regular pencil
\end{thmenumerate}
\end{lemma}

\newcommand{\Sop}{\mathsf{S}_{λ}}

\begin{pr}
\begin{enumerate}
\item
Consider, for any $λ\in\CC$, the operator $\Sop$ defined by
\[\Sop := λ\E + \A
.
\]
$\Sop$ can be decomposed into $\Sop'$ and $[\Sop]$ according to the following commuting diagram:
\begin{center}
\begin{tikzpicture}
\matrix(m) [commdiag]
{0 & \Mv' & \Mv & \Mv/\Mv' & 0 \\
0 & \Vv' & \Vv & \Vv/\Vv' & 0 \\};
\path[->]
(m-1-1) edge (m-1-2)
(m-1-2) edge (m-1-3)
(m-1-3) edge (m-1-4)
(m-1-4) edge (m-1-5)
(m-2-1) edge (m-2-2)
(m-2-2) edge (m-2-3)
(m-2-3) edge (m-2-4)
(m-2-4) edge (m-2-5)
(m-1-2) edge node[auto] {$\Sop'$} (m-2-2)
(m-1-3) edge node[auto] {$\Sop$} (m-2-3)
(m-1-4) edge node[auto] {$[\Sop]$} (m-2-4);
\end{tikzpicture}
\end{center}

Since both rows are exact sequence, out of the three operators $\Sop'$, $\Sop$ and $[\Sop]$, if two of them are invertible then the third one is. 
One easy way to prove this fact\footnote{This is a very general result that holds in other contexts as well, since one may also prove it by diagram chasing.} is by choosing bases in $\Mv$ and $\Vv$ which are compatible with the subspaces $\Mv'$ and $\Vv'$. The operator $\Sop$ is then represented by a block triangular matrix where the diagonal blocks are the matrices of $\Sop'$ and $[\Sop]$.
Now it is easy to check that if two of those three matrices are invertible, the third one is.

\item
Notice that for any $λ\in\CC$, $[\Sop] = [\A]$ (the operator $[\A]$ is defined in \autoref{prop_quotient_op}), so $[\Sop]$ is invertible if and only if $β_1^+ = 0$. 
As a result, we obtain the property
\[β^+_1 =0 \implies \bigl[ \forall λ\in\CC\quad \Sop \text{ invertible} \iff \Sop' \text{ invertible}\bigr]
.
\]
\item
For any $λ\in\CC$, the surjectivity of $\Sop$ implies that of $[\Sop]$. Since $[\Sop]$ does not depend on $λ$, it means that if $[\Sop]=[\A]$ is not surjective, then $\Sop$ is not surjective for any $λ\in\CC$. Now since, by definition, if $β^+_1≠0$ then $[A]$ is not surjective, we conclude that
\[β^+_1 ≠ 0 \implies \forall λ\in\CC\quad \Sop \text{ not invertible}
.
\]
\end{enumerate}
All the possibilities are covered and the claim is proved.
\end{pr}

\begin{pr}[Proof of \autoref{prop_reg_pencil}]
\begin{enumerate}
\item 
We first show by induction on the index that $(\E,\A)$ is a regular pencil if and only if $β^+=0$ and $(\E,\A)\red{∞}$ is a regular pencil. 
It is easy to show using \autoref{lma_regular_pencil}.
\item Now a system $(\E,\A)$ is a regular pencil if and only if the dual system $(\E,\A)^*$ is a regular pencil, so we may apply on $(\E\reds,\A\reds)$ the claim just proved.
Because of \autoref{def_betam}, we obtain that $(\E,\A)$ is a regular pencil if and only if $β^+$ and $β^-$ are zero, and $(\E\redss,\A\redss)$ is a regular pencil
\item Since, by \autoref{prop_Edyn_inv}, $\E\redss$ is invertible, the system $(\E\redss,\A\redss)$ is a regular pencil, and the claim is proved.
\end{enumerate}
\end{pr}

\section{Coupling}
\label{sec_coupling}

\subsection{Motivation: coupling spaces}
\label{coupling_spaces_sec}

In \autoref{sec_reduction} we showed how to define \emph{invariant} subspaces for the system $(\E,\A)$. ``Invariant'' means here that those subspaces are not arbitrarily chosen, they depend in a unique way from the system at hand.

In order to obtain a simple matrix representation of that system, we will need to choose supplementary spaces to the invariant subspaces $\Mv'$ and $\Vv'$.
In this section, we focus on such supplementary spaces for one reduction step only, and establish some results which will be needed in \autoref{sec_decomposition}.

We first look at the case of supplementary subspaces to the subspace $\Mv'$, i.e.,  subspaces $\Nv'\subset\Mv$ such that
\[\Mv = \Mv' \oplus \Nv'
.
\]

The strategy is to try and choose $\Nv'$ in the same direction as the part of $\ker\E$ that remains out of $\Mv'$.
First we define what this space is by decomposing the kernel of $\E$ in the part that is included in $\Mv'$ and some supplementary space.
This is achieved by choosing any supplementary space $\Kv'$ such that
\[\ker \E = (\ker\E \cap \Mv') \oplus \Kv'
.
\]
Then since, by construction, $\Kv' \cap \Mv' = 0$ one may complete $\Mv'$ by choosing a supplementary space {$\Cv'$} such that
\[\Mv = \Mv' \oplus \Cv' \oplus \Kv'
.
\]
We now define $\Nv'$ as
\[\Nv' := \Cv' \oplus \Kv'
.
\]

The choice of $\Cv'$ will prove to be essential to obtain a complete decomposition of the system $(\E,\A)$.
The tool to choose $\Cv'$ appropriately will be \autoref{lma_coupling_E}.

But notice now that no matter how we choose $\Cv'$, the space $\Kv'$ roughly speaking corresponds to the variables that are \emph{decoupled} from the rest of the system. 
They are sometimes called the \emph{algebraic constraints}.

\begin{example}
Let us illustrate the previous remark by a trivial example.
Consider the simple system
\[\begin{eqsys}
&x' = x\\
& y = 0
.
\end{eqsys}
\]
The variable $y$ is \emph{decoupled} from the rest of the system.
\end{example}



\subsection{Coupling Lemma for $\E$}


We will assume that some coupling space $\Wv''$ has already been chosen in the reduced system, and that will serve as a starting point for the choice of the coupling space at the present stage.
More precisely, we assume that the reduced space $\Vv'$ is already decomposed as
\[\Vv' = \E\Mv = \E\Mv' \oplus \Wv''
.
\]
That decomposition allows to construct the coupling spaces in an optimal manner, in one subspace $\Cv'$ coupled with $\Wv''$, and complement with vectors in the null-space of $\E$.

We state the result in a lemma, formulated outside the context of linear systems.

\begin{figure}
\begin{center}
\includegraphics[width=.6\textwidth]{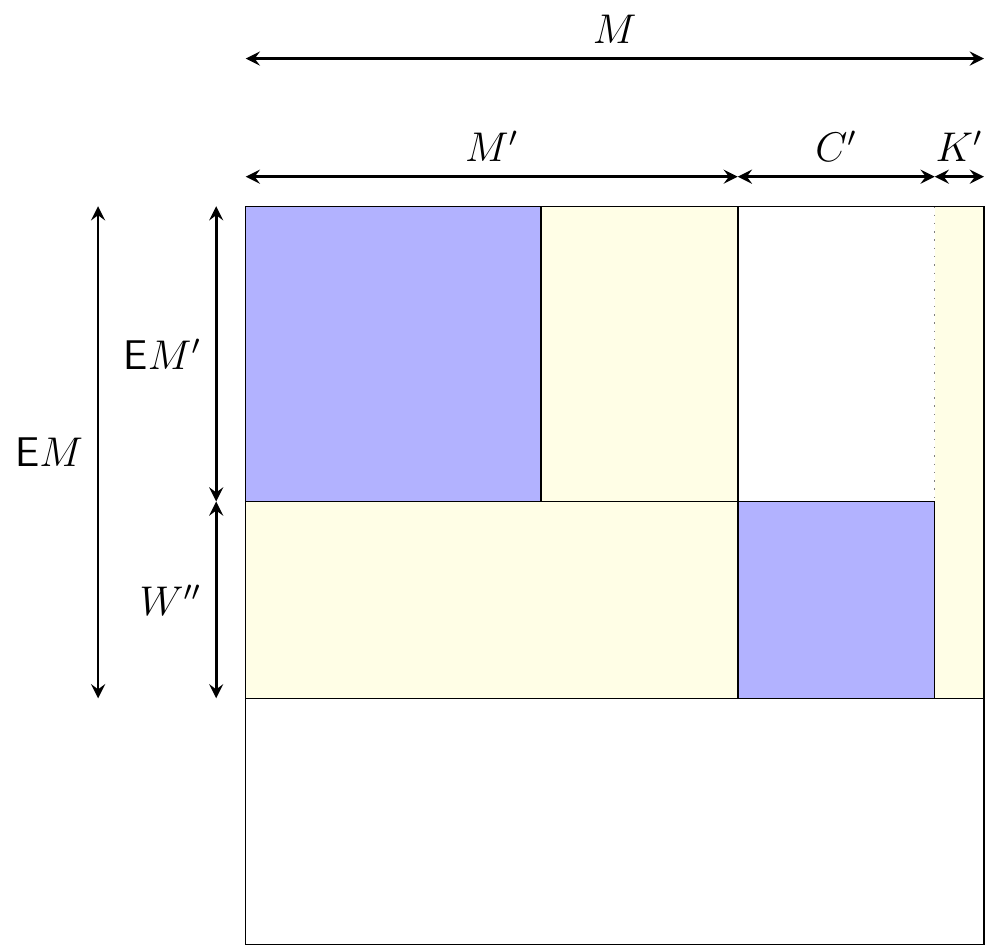}
\end{center}
\caption{An illustration of \autoref{lma_coupling_E}. There is a basis choice such that the operator $\E$ is represented as this matrix. Blue squares are identity matrix blocks. Other areas are filled with zeros.}
\label{fig_coupling_lemma}
\end{figure}

\begin{lemma}
\label{lma_coupling_E}
Assume that an operator $\E$ acting on a space $\Mv$, and consider a subspace $\Mv'\subset\Mv$. 
For any subspace $\Wv''$ such that $\E\Mv = \E\Mv' \oplus \Wv''$ there exists subspaces $\Kv'$ and $\Cv'$ such that
\[\Mv = \Mv' \oplus \Cv' \oplus \Kv'
\]
and such that the sequence
\begin{center}
\begin{tikzpicture}[scale=1,>=stealth]
\matrix(m) [exseq]
{0 & \Kv' & \Kv'\oplus\Cv' & \Wv'' & 0\\};
\path[->] (m-1-1) edge (m-1-2)
(m-1-2) edge (m-1-3)
(m-1-3) edge node[auto] {$\E$} (m-1-4)
(m-1-4) edge (m-1-5);
\end{tikzpicture}
\end{center}
is \emph{exact}. 
The exactness means here that $\ker\E\cap (\Kv'\oplus\Cv') = \Kv'$ and $\E(\Kv'\oplus\Cv') = \Wv''$.

Moreover, for any choice of basis in $\Wv''$ one may choose a basis of $\Cv'$ such that its image by $\E$ is the basis in $\Wv''$ (see \autoref{fig_coupling_lemma}).
\end{lemma}

\begin{pr}
\begin{enumerate}
\item
Consider
\[\overline{\Cv'} := \E\inv\Wv'' = \bigl\{x\in\Mv:\ \E x \in \Wv''\bigr\}
.
\]
Observe that $\Mv' + \overline{\Cv'} = \Mv$, and $\E\overline{\Cv'} = \Wv''$.
\item
Pick $x\in\overline{\Cv'}\cap\Mv'$. It implies that $\E x \in\Wv''\cap\E\Mv'$, so $\E x = 0$, i.e., $x\in\ker\E$. 
We conclude that $\overline{\Cv'}\cap\Mv' \subset \ker\E$.
\item Choose $\Cv'$ such that $\overline{\Cv'} = \ker\E \oplus \Cv'$. 
It follows from the previous observation that $\Cv'\cap \Mv' = 0$. This implies
\[\Mv = (\Mv' + \ker\E) \oplus \Cv'
.
\]
\item Decompose further $\ker\E$ as
\[\ker\E = (\ker\E\cap\Mv')\oplus \Kv'
.
\]
As a consequence, we obtain
\[\Mv' + \ker\E = \Mv' \oplus \Kv'
\]
and $\Kv' \subset \ker\E$.
\item
Finally we have $\E\Cv' = \E\overline{\Cv'} = \Wv''$, and $\Cv'\cap\ker\E$ = 0. 
Moreover, since $\E$ restricted on $\Cv'$ sends $\Cv'$ bijectively to $\Wv'$, the inverse image of the basis of $\Wv''$ is a basis of $\Cv'$.
\end{enumerate}

\end{pr}

\subsection{Complementary subspaces}

The definition of $\Mv'$, in \autoref{space_red_def}, is the set of vectors in $\Mv$ such that $\A x$ intersects with the image of $\E$.
We are now interested in the converse statement, namely, that if we choose a subspace $\Nv'$ which does not intersect $\Mv'$, its image $\A\Nv'$ by $\A$ should not intersect the image of $\E$.
This is the gist of the following lemma.

\begin{lemma}
\label{lma_AN_imE}
Consider a linear system $(\E,\A)$ and assume that $\Nv'$ is a subspace of $\Mv\EA$ that does not intersect $\Mv'\EA$, i.e.
such that
\[\Nv' \cap \Mv'\EA = 0
.
\]
Then the property
\[\A \Nv'\cap \Vv'\EA = 0
\]
holds.
\end{lemma}

Note that this result is a consequence of \autoref{prop_quotient_op}.
Indeed, the subspace $\Nv'$ may be injected in $\Mv/\Mv'$ and the result follows from the fact that $[\A]$ is injective.

We give also a direct proof of this elementary lemma.

\begin{pr}
For general subspaces $\Nv'\subset\Mv$ and $\Wv'\subset\Vv$, we have
\[\A\Nv' \cap \Wv' = \A(\Nv' \cap \A\inv\Wv')
.
\]
With $\Wv' = \Vv'$ and since by \autoref{space_red_def}, $\Mv' = \A\inv\Vv'$, we obtain
\[\A\Nv' \cap \Vv' = \A(\Nv' \cap \Mv')
\]
from which the claim follows.
\end{pr}

\subsection{Coupling Lemma for Systems}

We may now combine \autoref{lma_coupling_E} and \autoref{lma_AN_imE} and obtain a fundamental Lemma that decouples the operators $\E$ and $\A$ on supplementary spaces to $\Mv'$ and $\Vv'$.
 
\begin{lemma}
\label{lma_coupling_system}
Suppose that there is a decomposition
\[
\Vv' = \Vv'' \oplus \Wv''
,
\]
and that $\Wv''$ is equipped with a basis.

Then there exists decompositions
\begin{align*}
\Mv &= \Mv' \oplus \Nv'
, \\
\Vv &= \Vv' \oplus \Wv'
,
\end{align*}
and subspaces
\[\begin{aligned}
\Cv'&\subset\Mv \qquad &\Dv'&\subset\Vv
,\\
\Kv'&\subset\Mv \qquad &\Zv'&\subset\Vv,
\end{aligned}\]
 such that
\begin{align}
\Nv' &= \Cv' \oplus \Kv'
,\\
\Wv' &= \Dv' \oplus \Zv'
.
\end{align}

Those subspaces are such that the following sequences are exact:

\begin{center}
\begin{tikzpicture}
\matrix(m) [exseq]
{0 & \Kv'  & \underbrace{\Kv' \oplus \Cv'}_{\Nv'} & \Wv'' & 0\\
0 & \Nv' & \underbrace{\Dv' \oplus \Zv'}_{\Wv'} & \Zv' & 0\\};
\path[->]
(m-1-1) edge (m-1-2)
(m-1-2) edge (m-1-3)
(m-1-3) edge node[auto] {$\E$} (m-1-4)
(m-1-4) edge (m-1-5);
\path[->]
(m-2-1) edge (m-2-2)
(m-2-2) edge node[auto] {$\A$} (m-2-3)
(m-2-3) edge (m-2-4)
(m-2-4) edge (m-2-5);
\end{tikzpicture}
\end{center}
and such that
\[\A\Mv \cap \Zv' = 0
.
\]

Moreover, one may choose basis in the subspaces $\Cv'$, $\Kv'$, $\Dv'$ and $\Zv'$ such that the basis of $\Dv'$ is the image by $\A$ of the basis of $\Nv'$, and the basis on $\Wv''$ is the image by $\E$ of the basis of $\Cv'$.

\end{lemma}
\begin{pr}
\begin{enumerate}
\item
By the assumption on $\Wv''$, we have
\[\E\Mv = \Vv' = \Vv'' \oplus \Wv'' = \E\Mv' \oplus \Wv''
.
\]
The subspace $\Wv''$ is moreover equipped with a basis by the induction hypothesis.
\item Appealing to \autoref{lma_coupling_E} we obtain subspaces $\Cv'$ and $\Kv'$ such that
\[\Mv = \Mv' \oplus \Cv' \oplus \Kv'
\]
with
\[\E \Cv' = \Wv''
\]
and
\[\E \Kv' = 0
\]
and
\[\ker \E \cap \Cv' = 0
.
\]
Note that, given a basis in $\Wv''$ we can choose a basis on $\Cv'$ such that $\E$ sends that basis on that of $\Wv''$.

Let us now define the subspace $\Nv'\subset\Mv$ by
\[\Nv' := \Cv' \oplus \Kv'
.
\]
We choose an arbitrary basis of the space $\Kv'$, and this provides us with a basis for the space $\Nv'$.
\item Recall now that, according to \autoref{lma_AN_imE},
\[\A \Nv' \cap \E\Mv = 0
,
\]
and by \autoref{kerA_sub_Mprim_prop}, the operator $\A$ sends the basis of the space $\Nv'$ to a set of independent vectors in the space $\Vv$.
We thus choose as a basis of $\A \Nv'$ the image of the basis of $\Nv'$ by $\A$.
\item 
Now choose a subspace $\Zv'\subset\Vv$ such that
\[\Vv = \E\Mv \oplus \A \Nv' \oplus \Zv'
\]
and pick an arbitrary basis of that subspace.
We define
\[\Dv' := \A\Nv'\]
and
\[\Wv' := \Dv' \oplus \Zv'
.
\]
\todo{check the claims of the theorem?}
\end{enumerate}

\end{pr}

\begin{remark}
\label{prop_dimKZ_ab}
The dimensions of the spaces introduced in \autoref{thm_linear_decomposition} are related to the dimensions of the spaces introduced in \autoref{def_DMDV}, and to the defects (\autoref{def_defects}).
The relations are given by
\[\begin{aligned}
\dim\Wv'&= \dim\DV' \qquad &\dim\Zv'&= α_1
,\\
\dim\Nv'&= \dim\DM' \qquad &\dim\Kv'&= β^+_1
.
\end{aligned}\]
\end{remark}

\section{Strangeness}
\label{sec_strangeness}

In order to illustrate the power of reduction, and to show an application of \autoref{lma_coupling_system}, we show an intermediate result.
Instead of looking at the equivalence classes for the equivalence of matrices, that is, pairs of invertible operators acting on $(\E,\A)$ as $(\Pm\E\Qm,\Pm\A\Qm)$, we look at the \emph{weak} equivalence.

Weak equivalence is determined by another group, which elements consist of two invertible operators $\Pm$ and $\Qm$ and an arbitrary operator $\mat{R}$, acting on a system $(\E,\A)$ as
\begin{equation}
\label{weak_equiv_action_eq}
(\Pm,\Qm,\mat{R}) \cdot (\E,\A) := (\Pm\E\Qm,\Pm(\E\Rm + \A\Qm))
.
\end{equation}

\subsection{Weak equivalence group}

The group operation corresponding to weak equivalence is given by
\[(\Pm_2,Q_2,R_2) \cdot (\Pm_1,Q_1,R_1) = (\Pm_2\Pm_1,\Qm_1 \Qm_2,Q_1R_2+R_1Q_2)
,
\]
where $\Pm_1$, $\Pm_2$ are automorphisms of $\Vv$, $\Qm_1$, $\Qm_2$ are automorphisms on $\Mv$, and $\Rm_1$, $\Rm_2$ are arbitrary endomorphisms on $\Mv$.

The identity is then
\[(I,I,0)
,
\]
and the inverse of an element $(\Pm,\Qm,\Rm)$ is given by
\[(\Pm,\Qm,\Rm)\inv = (\Pm\inv,Q\inv,-Q\inv R Q\inv)
.
\]

Clearly, the elements of the form $(\Pm,\Qm,0)$ form a subgroup corresponding to the equivalence relation.
Another subgroup is given by elements of the form $(\Id,\Id,\Rm)$.

For the study of the orbits of the weak equivalence group, the identity
\begin{equation}
\label{orbit_sep_eq}
(\Pm,Q,R) = (\Pm,Q,0)\cdot (I,I,R Q\inv) = (I,I,Q\inv R) \cdot (\Pm,Q,0)
\end{equation}
 shows that we may restrict our attention to one subgroup at a time.

\subsubsection{Orbit Invariants}

The orbits of the weak equivalence group action \eqref{weak_equiv_action_eq} were studied in \cite{Mehrmann}, in which the authors exhibited a complete set of invariants.
We give an alternative proof here, thereby shedding some light on the notion of \emph{strangeness}.

\begin{theorem}
\label{thm_strangeness}
	A complete set of invariants for the group action \eqref{weak_equiv_action_eq} is given by
	\begin{enumerate}
		\item $d := \dim \Vv''$,
		\item $a := \dim \ker [\E] = α_1$,
		\item $s := \dim \DV''$.
	\end{enumerate}
\end{theorem}

The integer $s$ is called ``strangeness'' in \cite{Mehrmann}.

\begin{figure}
\begin{center}
\includegraphics[width=.5\textwidth]{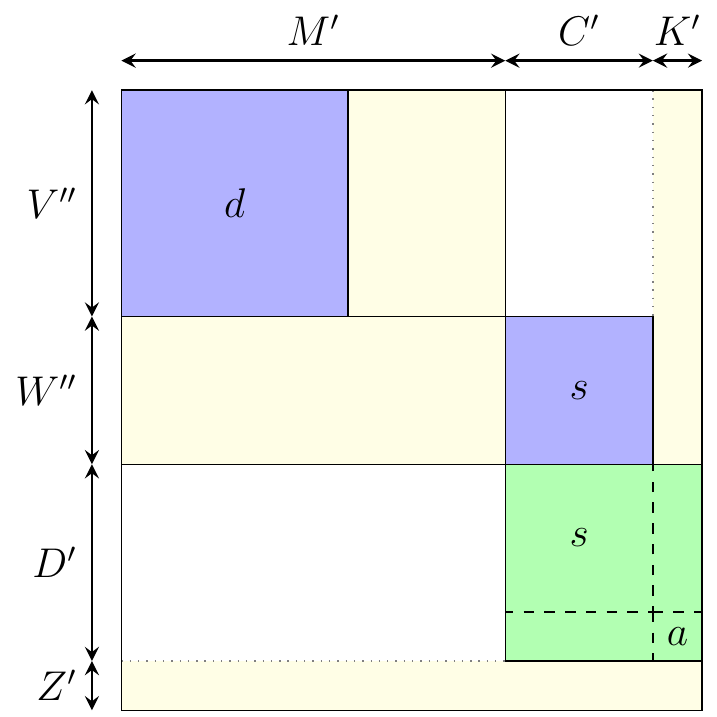}
\end{center}
\caption{Canonical form of a matrix corresponding to the weak equivalence. The matrix $\E$ is represented in blue, whereas the matrix for $\A$ is represented in green. All such squares are identity matrices. The rest is filled out by zero entries.}
\label{fig_strangeness}
\end{figure}

\begin{pr}
\begin{enumerate}
\item
First we have to check that the three integers are indeed invariants of the group action.
Clearly, they are invariants by transformations of the form $(\Pm,\Qm,0)$, which are merely equivalent transformation.

Let us examine the case of a transformation
\[
(\Eb,\Ab) = (\Id,\Id,\Rm)\cdot (\E,\A) = (\E, \E R + \A)
.
\]
We have
\[
\overline{\Vv}' = \Eb \Mv = \E \Mv = \Vv'
,
\]
\[
\overline{\Mv}' := \{x:\ \Ab x \in\Eb \Mv\} = \Mv'
,
\]
so
\[
(\Eb',\Ab') = (\E,\A)
\]
and
\[
\Vv'' = \Eb'\Mv' = \E\Mv' = \Vv''
.
\]
Using \eqref{orbit_sep_eq}, this shows that the spaces $\Vv'$, $\Vv''$ and $\Mv'$ are invariants of all of the weak equivalence group transformations.

As a result, the spaces $∆\Vv''$ and the operator $[\E]$ are also invariants, so the integers $d$, $a$ and $s$ are invariants.
\item
Now we show that the integers $d$, $a$ and $s$ are the only invariants.
In order to show that, we show that a system $(\E,\A)$ is weakly equivalent to a canonical form that depends only on those three integers.

In order to achieve this, we decompose $\Mv$ and $\Vv$ using \autoref{lma_coupling_system}.
\begin{enumerate}
\item
Let us choose an arbitrary decomposition
\[
\Vv' = \Vv'' \oplus \Wv''
.
\]
We may now apply \autoref{lma_coupling_system} to obtain spaces $\Cv'$, $\Kv'$, $\Dv'$ and $\Zv'$ equipped with appropriate bases.
Using \autoref{prop_dimKZ_ab} we obtain $s = \dim \Wv''$ and $a = \dim \Kv'$.
\item
Finally, define $\Pi$ as a projector from $\Mv$ to $\Mv'$ along $\Nv'$. 
Let $\mat{F}$ be a right inverse for $\E$ on $\Vv' = \E\Mv$. 

Define
\[\Rm := -\mat{F}\A \Pi
,
\]
so $\E\Rm + \A = 0$ on $\Mv'$.

As a result, if we define the new system $(\Eb,\Ab)$ by
\[(\Eb,\Ab) := (\E,\E\Rm+\A)
,
\]
then the restriction of $\Ab$ on $\Mv'$ is zero.
\item
Now we may choose a basis of $\Mv'$ and of $\Vv''=\E\Mv'$ such that $\Eb$ is represented by the identity matrix on $\Mv'$. 

This provides us with complete basis of $\Mv$ and $\Vv$ such that the matrices $\Eb$ and $\Ab$ take the form described in \autoref{fig_strangeness}.
\end{enumerate}
\end{enumerate}
\end{pr}

\section{Direct Decomposition}
\label{sec_decomposition}

\subsection{Decomposition Theorem}

In \autoref{sec_reduction} we showed how to define invariant subspaces $\Mv\red{k}$ and $\Vv\red{k}$ for the system $(\E,\A)$.
In order to obtain a complete decomposition of the spaces $\Mv$ and $\Vv$, it is necessary to construct subspaces that bridge the gap between each invariant subspaces $\Mv\red{k}$ and $\Vv\red{k}$. 
More precisely we construct spaces $\Nv\red{k}$ and $\Wv\red{k}$ such that
\[\Mv\red{k} = \Mv\red{k+1} \oplus \Nv\red{k+1}
\text{,} \qquad 
\Vv\red{k} = \Vv\red{k+1} \oplus \Wv\red{k+1}
.
\]
In a sense, the subspaces $\Nv\red{k}$ and $\Wv\red{k}$ correspond to the spaces $\DM\red{k}$ and $\DV\red{k}$ respectively, defined in \autoref{def_DMDV}.

The construction of those supplementary spaces proceeds backwards, in the direction opposite to the reduction.
One must first totally reduce the system $(\E,\A)$.
Assume that the index is $n$.
One then chooses an arbitrary complementary space $\Wv\red{n}$ such that $\Vv\red{n-1} = \Vv\red{n}\oplus\Wv\red{n}$, and equip that space with an arbitrary basis.
The rest is a repeated application of \autoref{lma_coupling_system}.

We now see how to choose those supplementary spaces $\Nv\red{k}$ and $\Wv\red{k}$ so that the operators $\E$ and $\A$ are \emph{simultaneously} decomposed in an advantageous way.

\begin{theorem}
\label{thm_linear_decomposition}
Consider a system $(\E,\A)$.

Recall the definitions of the subspaces $\Mv\red{k}$ and $\Vv\red{k}$ in \eqref{space_red_notation}.

For any integer $k \in \NN$ there exists  subspaces $\Nv\red{k+1}\subset\Mv$ and $\Wv\red{k+1}\subset\Vv$ such that

\begin{align*}
\Mv\red{k} &= \Mv\red{k+1} \oplus \Nv\red{k+1}
, \\
\Vv\red{k} &= \Vv\red{k+1} \oplus \Wv\red{k+1}
,
\end{align*}

and for any integer $k≥1$ there exists subspaces
\[\begin{aligned}
\Cv\red{k}&\subset\Mv \qquad &\Dv\red{k}&\subset\Vv
,\\
\Kv\red{k}&\subset\Mv \qquad &\Zv\red{k}&\subset\Vv,
\end{aligned}\]
 such that
\begin{align}
\Nv\red{k} &= \Cv\red{k} \oplus \Kv\red{k}
,\\
\Wv\red{k} &= \Dv\red{k} \oplus \Zv\red{k}
.
\end{align}

Those subspaces are such that for any integer $k≥1$, the following sequences are exact (see \autoref{fig_decomposition}):

\begin{center}
\begin{tikzpicture}
\matrix(m) [exseq]
{0 & \Kv\red{k}  & \underbrace{\Kv\red{k} \oplus \Cv\red{k}}_{\Nv\red{k}} & \Wv\red{k+1} & 0\\
0 & \Nv\red{k} & \underbrace{\Dv\red{k} \oplus \Zv\red{k}}_{\Wv\red{k}} & \Zv\red{k} & 0\\};
\path[->]
(m-1-1) edge (m-1-2)
(m-1-2) edge (m-1-3)
(m-1-3) edge node[auto] {$\E$} (m-1-4)
(m-1-4) edge (m-1-5);
\path[->]
(m-2-1) edge (m-2-2)
(m-2-2) edge node[auto] {$\A$} (m-2-3)
(m-2-3) edge (m-2-4)
(m-2-4) edge (m-2-5);
\end{tikzpicture}
\end{center}
and such that
\[\A\Mv\red{k-1} \cap \Zv\red{k} = 0
.
\]

Moreover, one may choose basis in the spaces $\Cv\red{k}$, $\Kv\red{k}$, $\Dv\red{k}$ and $\Zv\red{k}$ such that the basis of $\Dv\red{k}$ is the image by $\A$ of the basis of $\Nv\red{k}$, and the basis on $\Wv\red{k+1}$ is the image by $\E$ of the basis of $\Cv\red{k}$.

\end{theorem}

\begin{remark}
For the reader averse to the language of exact sequences, the fact that the sequences of \autoref{thm_linear_decomposition} are \emph{exact} means in that case that
\begin{equation*}
\begin{aligned}
\E\Cv\red{k} &= \Wv\red{k+1}
,\\
\E\Kv\red{k} &= 0
,\\
\A\Nv\red{k} &= \Dv\red{k}
,\\
\ker\A\cap\Nv\red{k}&=0
.\\
\end{aligned}
\end{equation*}

\end{remark}

\begin{remark}
\label{prop_dimKZ_ab_k}
In the same spirit as \autoref{prop_dimKZ_ab}, we notice the relation between the dimensions of the various subspaces introduced in \autoref{thm_linear_decomposition}, and the dimensions of the spaces defined in \autoref{def_DMDV}, and to the defects (\autoref{def_defects}).
For any integer $k≥1$, the relations are given by
\[\begin{aligned}
\dim\Wv\red{k}&= \dim\DV\red{k} \qquad &\dim\Zv\red{k}&= α_k
,\\
\dim\Nv\red{k}&= \dim\DM\red{k} \qquad &\dim\Kv\red{k}&= β^+_k
.
\end{aligned}\]
\end{remark}

\matrixfigure[.5]{fig_decomposition}{decomposition}{An illustration of the decomposition described in \autoref{thm_linear_decomposition}. 
Noticing that no matter what bases we choose in $\Mv\red{∞}$ and $\Vv\red{∞}$ the matrix is block diagonal (see \autoref{cor_operator_splitting}), we may choose those bases in such a way that $\E$ is represented as the identity matrix on that block. 
Since we now that $\E\red{∞}$ is surjective (\autoref{prop_Eredinf_surj}), that identity block stretches to fill $\Vv\red{∞}$.}

\matrixfigure{fig_induction_proof}{proof}{An illustration of \autoref{thm_linear_decomposition} and \autoref{lma_coupling_system} on an index three system.
The grey shaded part pictures the previous step of the recursion.
Starting with $\Wv''$, one constructs the space $\Cv'$ such that $\E\Cv'=\Wv''$, and a subspace $\Kv'$ such that $\Kv'\subset\ker\E$ using \autoref{lma_coupling_E}, and defines $\Nv' := \Cv'\oplus \Kv'$.
One then constructs $\Zv'$ such that $\Vv = \Vv' \oplus \A \Nv' \oplus \Zv'$.
This in turn defines $\Wv':=\A\Nv'\oplus\Zv'$.}

\begin{pr}[Proof of \autoref{thm_linear_decomposition}]
We  proceed by induction on the index (see \autoref{fig_induction_proof}).
If the index is zero, all the spaces $\Mv\red{k}$ and $\Vv\red{k}$ are zero, and there is nothing to prove.
Assume now that the statement holds for systems of index $n-1$.
Given a system $(\E,\A)$ of index $n$, the reduced system $(\E',\A')$ has index $n-1$, so we may apply the induction hypothesis on that reduced system.

For clarity, let us denote
\[(\overline{\E},\overline{\A}):= (\E',\A')
,\qquad
\overline{\Mv} = \Mv'
,\quad 
\overline{V}=\Vv'
.
\]
The reduced system $(\E,\A)'$ consists of operators operating from $\Mv'$ to $\Vv'$, so by the induction hypothesis we obtain a decomposition of the spaces $\overline{\Mv}$ and $\overline{\Vv}$ into subspaces $\overline{\Wv}\red{k}$, $\overline{\Nv}\red{k}$ as described in the statement of the theorem.

We have to shift the indices of all the spaces produced for the final statement to hold.
For example, we define for any integer $k≥2$
\[
\Wv\red{k} := \overline{\Wv}\red{k-1}
,
\]
so we may write the decomposition of $\Vv'$ as
\[\overline{\Vv} = \Vv' = \Vv\red{∞} \oplus \Wv\red{n} \oplus \Wv\red{n-1}\oplus \cdots  \oplus \Wv''
.
\]

The reduced operators $\E'$ and $\A'$ being restrictions of $\E$ and $\A$, the statements obtained from the induction hypothesis apply to the operators $\E$ and $\A$.
\todo{clearer explanation}

Applying \autoref{lma_coupling_system} yields the desired result.

\end{pr}

\subsection{Decomposition in invariant subspaces}
\todo{change name to invariant subsystems?}

A crucial consequence of \autoref{thm_linear_decomposition} is that it provides us with decompositions of $\Mv$ and $\Vv$ such that $\E$ and $\A$ may be restricted on those subspaces:
\begin{corollary}
\label{cor_operator_splitting}
\todo{where is this corollary used?}
Given the decomposition provided by \autoref{thm_linear_decomposition}, and defining $\overline{\Mv}$ and $\overline{\Vv}$ by
\[
	\overline{\Mv} := \bigoplus_{k}\Nv\red{k}
	, \qquad
	\overline{\Vv} := \bigoplus_{k}\Wv\red{k}
	,
\]
then, by construction,
\[
	\Mv = \Mv\red{∞} \oplus \overline{\Mv}
	, \qquad
	\Vv = \Vv\red{∞} \oplus \overline{\Vv}
	,
\]
and we have
\begin{thmenumerate}
	\item \[\E\Mv\red{∞} = \Vv\red{∞}\qquad \A\Mv\red{∞}\subset\Vv\red{∞}\]
	\item \[\E\overline{\Mv}\subset\overline{\Vv}\qquad\A\overline{\Mv}\subset\overline{\Vv}\]
\end{thmenumerate}
\end{corollary}
\todo{simple illustration of the splitting!!}
\begin{pr}The fact that $\E\overline{\Mv}\subset\overline{\Vv}$ and $\A\overline{\Mv}\subset\overline{\Vv}$ follows from
\[\E\Nv\red{k} \subset \Wv\red{k+1} \subset\overline{\Vv}
\]
and
\[\A\Nv\red{k}\subset\Wv\red{k} \subset\overline{\Vv}
.
\]
\end{pr}

\section{Dual Decomposition}
\label{sec_decomposition_bis}

\todo{mention transpose of a matrix}

\subsection{Dual space decomposition}

Assume that a finite dimensional vector space $\Mv$ is decomposed in a direct sum of subspaces, i.e.,
\[\Mv = \Mv_1 \oplus \Mv_2 \oplus \cdots \oplus \Mv_n
.
\]

This decomposition induces the dual space decomposition
\[(\Mv_k)_* := \bigl( \bigoplus_{j≠k} \Mv_j \bigr)^{\perp} = \Bigl\{ φ\in\Mv^*:\ \bracket{φ}{x} = 0\ \forall x\in \bigoplus_{j≠k} \Mv_j \Bigr\}
.
\]

Although it is not reflected by the notation, it is clear that $(\Mv_k)_*$ actually depends not only on $\Mv_k$ but on all the other spaces of the decomposition.
It is a generalization of the notion of dual basis.

Assume further that $\Mv$ is equipped with a basis.
We say that this basis is \emph{compatible} with the decomposition if each subspace is the span of a subset of the basis.

If $\Mv$ is equipped with a basis $\Bb$ compatible with a subspace decomposition, then the dual basis is compatible with the dual space decomposition.

Indeed, consider a subspace $\Mv_k$ of the decomposition. 
Since the basis is compatible with the decomposition, that subspace is spanned by a subset of the basis $\Bb$, say $\Sb_{\Mv_k}\subset\Bb$, i.e.,
\[\Mv_k = \Span \Sb_{\Mv_k}
.
\]
The dual decomposition is such that the associated subspace $(\Mv_k)_*$ is the span of the dual basis with the \emph{same} subset $\Sb_{\Mv_k}$, i.e.,
\[(\Mv_k)_* = \Span \bigl\{e^*: e\in \Sb_{\Mv_k}\bigr\}
.
\]
Here the covector $e^*$ is the element of the dual basis of $\Bb$ corresponding to $e$, i.e., such that
\[\bracket{e^*}{f} = \begin{cases}
	0 & \text{if $f≠e$}\\
	1 & \text{if $f = e$}
\end{cases}
\]

\begin{lemma}
\label{lma_exseq_flip}
Assume that $A$ and $B$ are subspaces of $\Mv$ that are part of a subspace decomposition of $\Mv$, and that $C$ is a subspace of $\Vv$ that is part of a subspace decomposition of $\Vv$. Assume further that $\Mv$ and $\Vv$ are equipped with bases compatible with their decompositions. Take an operator $\mat{S}$ operating from $\Mv$ to $\Vv$.

The following two statements are equivalent:
\begin{thmenumerate}
\item
The following sequence is exact:
\begin{center}
\begin{tikzpicture}
\matrix(m) [exseq]
{0 & A  & A\oplus B & C & 0\\};
\path[->]
(m-1-1) edge (m-1-2)
(m-1-2) edge  (m-1-3)
(m-1-3) edge node[auto] {$\mat{S}$} (m-1-4)
(m-1-4) edge (m-1-5);
\end{tikzpicture}
\end{center}
Moreover, the operator $\mat{S}$ sends the basis of $B$ on the basis of $C$.

\item
The following ``dual'' sequence is exact:
\begin{center}
\begin{tikzpicture}
\matrix(m) [exseq]
{0 & C_*  & A_*\oplus B_* & A_* & 0\\};
\path[->]
(m-1-1) edge (m-1-2)
(m-1-2) edge node[auto] {$\mat{S}^*$} (m-1-3)
(m-1-3) edge  (m-1-4)
(m-1-4) edge (m-1-5);
\end{tikzpicture}
\end{center}
and
\[\mat{S}^*\Vv^*\cap \A_* = 0
.
\]
Moreover, the operator $\mat{S}^*$ sends the basis of $C_*$ on the basis of $B_*$.
\end{thmenumerate}
\end{lemma}
\begin{pr}
Denote the basis on $\Mv$ and $\Vv$ by $\Bb(\Mv)$ and $\Bb(\Vv)$ respectively. 
Clearly, for any $e\in\Bb(\Mv)$ and $f\in\Bb(\Vv)$, we have
\[ \bracket{f^*}{\mat{S} e} = \bracket{\mat{S}^* f^*}{e} = \bracket{(e^*)^*}{\mat{S}^* f^*}
\]
where $(e^*)^*$ is the dual basis of the dual basis of $\Bb$. 

The proof is now a simple verification by expressing each of the statements in terms of the bases. For example, $\mat{S}A = 0$ may be written as
\[\bracket{f^*}{\mat{S}e} = 0 \qquad \forall e\in\Sb_{A}\quad f\in\Bb(\Vv)
,
\]
so one obtains
\[\bracket{(e^*)^*}{\mat{S}^* f^*}=0 \qquad \forall e\in\Sb_{A}\quad f\in\Bb(\Vv)
,
\]
which means that $\mat{S}A = 0 \iff \mat{S}^* \Vv^* \cap A_* = 0$.
The other statements are verified in the same fashion.
\end{pr}

\subsection{Conjugate Decomposition}

\newcommand{\redd}[1]{\red{#1}_*}

Consider a finite dimensional linear system $(\E,\A)$ and its dual $(\E^*,\A^*)$. 
The corresponding domain and codomain are denoted by
\[\begin{split}
	\overline{\Mv} := \Mv\cdom{^*} = \Vv^*%
,\\
	\overline{\Vv} := \Vv\cdom{^*} = \Mv^*%
.
\end{split}\]

By applying \autoref{thm_linear_decomposition} on the dual system $(\E^*,\A^*)$ one obtains a decomposition of $\overline{\Mv}$ and $\overline{\Vv}$ as
\[\overline{\Mv} = \overline{\Mv}\red{∞} \oplus \overline{\Kv}' \oplus \overline{\Cv}' \oplus \cdots
,
\]
\[\overline{\Vv} = \overline{\Vv}\red{∞} \oplus \overline{\Zv}' \oplus \overline{\Dv}' \oplus \cdots
.
\]

Moreover, all those subspaces are equipped with a suitable basis.
By choosing a basis for the spaces $\overline{\Mv}\red{∞}$ and $\overline{\Vv}\red{∞}$ we obtain compatible bases $\Bb(\overline{\Mv})$ and $\Bb(\overline{\Vv})$ of $\overline{\Mv}$ and $\overline{\Vv}$ respectively.

\begin{theorem}
\label{thm_linear_decomposition_conj}
Consider the decompositions produced by \autoref{thm_linear_decomposition} for the dual system $(\E,\A)^*$.
The dual decompositions induce decompositions of the spaces $\Mv$ and $\Vv$ by the canonical isomorphism between a space and its bidual. 

For any integer $k \in \NN$ we have
\begin{align*}
\Mv\redd{k} &= \Mv\redd{k+1} \oplus \Nv\redd{k+1}  \\
\Vv\redd{k} &= \Vv\redd{k+1} \oplus \Wv\redd{k+1}  ,
\end{align*}


Those subspaces are such that the following sequences are exact (see \autoref{fig_decomposition}):

\todo{illustration}
\begin{center}
\begin{tikzpicture}
\matrix(m) [exseq]
{0 & \Wv\redd{k+1}  & \underbrace{\Kv\redd{k} \oplus \Cv\redd{k}}_{\Nv\redd{k}} & \Kv\redd{k} & 0\\
0 & \Zv\redd{k} & \underbrace{\Dv\redd{k} \oplus \Zv\redd{k}}_{\Wv\redd{k}} & \Nv\redd{k} & 0\\};
\path[->]
(m-1-1) edge (m-1-2)
(m-1-2) edge node[auto] {$\E$} (m-1-3)
(m-1-3) edge  (m-1-4)
(m-1-4) edge (m-1-5);
\path[->]
(m-2-1) edge (m-2-2)
(m-2-2) edge  (m-2-3)
(m-2-3) edge node[auto] {$\A$} (m-2-4)
(m-2-4) edge (m-2-5);
\end{tikzpicture}
\end{center}

Moreover, the basis of $\Cv\redd{k}$, $\Kv\redd{k}$, $\Dv\redd{k}$ and $\Zv\redd{k}$  are such that the basis of $\Nv\redd{k}$ is the image by $\A$ of the basis of $\Dv\redd{k}$, and the basis of $\Cv\redd{k}$ is the image by $\E$ of the basis of $\Wv\redd{k+1}$.

\end{theorem}

\begin{pr}
It is a direct application of \autoref{lma_exseq_flip}.
\end{pr}

\begin{remark}
	The exact sequences of \autoref{thm_linear_decomposition} are the same as those of \autoref{thm_linear_decomposition_conj} but with flipped arrows.
	It just reflects how the block structure of a matrix is related to the block structure of the transposed matrix.
\end{remark}

\subsection{Second sweep of the decomposition}

Remember that the constraint defects of the adjoint of the totally reduced system $(\E\reds,\A\reds)$  are zero (\autoref{prop_alpha_reds_zero}). Along with \autoref{prop_dimKZ_ab_k}, we conclude that the corresponding subspaces $\Kv\red{k}$ produced in \autoref{thm_linear_decomposition} are zero.


Now, using \autoref{cor_operator_splitting}\todo{really?}, we are in a position to use the decomposition of \autoref{thm_linear_decomposition} for the dual system $(\E\reds,\A\reds)$ and obtain a decomposition of $\Mv\reds$ and $\Vv\reds$.

\begin{theorem}
\label{thm_full_linear_decomposition}
In addition to the decomposition given by \autoref{thm_linear_decomposition}, the spaces $\Mv\red{∞}$ and $\Vv\red{∞}$ may now be decomposed as (see \autoref{fig_full_decomposition}):
\begin{center}
\begin{tikzpicture}
\matrix(m) [exseq]
{0 & \Zv\red{k} & \underbrace{\Zv\red{k} \oplus \Dv\red{k}}_{\Wv\red{k}} & \Nv\red{k} & 0 \\
0 & \Wv\red{k+1} &  \Zv\red{k} & 0& \\};
\path[->]
(m-1-1) edge (m-1-2)
(m-1-2) edge (m-1-3)
(m-1-3) edge node[auto] {$\A$} (m-1-4)
(m-1-4) edge (m-1-5);
\path[->]
(m-2-1) edge (m-2-2)
(m-2-2) edge node[auto] {$\E$} (m-2-3)
(m-2-3) edge (m-2-4);
\end{tikzpicture}
\end{center}

\end{theorem}

\begin{pr}

\end{pr}

\matrixfigure[.6]{fig_full_decomposition}{full}{
An illustration of the full decomposition of \autoref{thm_full_linear_decomposition}.
The first decomposition leads to $\Mv''$ and the corresponding space $\Vv''' = \E\Mv''$, at which point the  algorithm stalls.
The second step consists in transposing the reduced operators $\E\red{∞}$ and $\A\red{∞}$ (indicated by the bold red frame on the figure), running the same algorithm, and transposing back again.
The upper left checkered area denotes the identity for $\E$, and a non specific matrix for $\A$.
Notice that this block is completely separated from the rest, so one may now reduce the $\A$ to Jordan blocks by a similarity transformation.}

\begin{remark}
The various defects defined in \autoref{def_defects} and \autoref{def_betam} may now be pictured clearly using the \autoref{thm_full_linear_decomposition}; see \autoref{fig_kronecker}.
\end{remark}

\todo{add more typical examples: DAE, control, observation}


\section{Kronecker Indices}
\label{sec_kronecker}

\subsection{Basis Arrangement}

In this section we prove a result on the basis obtained in \autoref{thm_linear_decomposition}, which will be useful to determine  the relation with the Kronecker decomposition theorem.

\newcommand{\mv}{\mat{m}}
\newcommand{\vv}{\mat{v}}

\newcommand{\Eedge}[2]{\path[<-] (m-1-#1) edge node[auto] {$\E$} (m-1-#2);}
\newcommand{\Aedge}[2]{\path[->] (m-1-#1) edge node[auto] {$\A$} (m-1-#2);}
\newcommand{\EAedge}[3]{\Eedge{#1}{#2} \Aedge{#2}{#3}}

\begin{definition}
\label{def_NL_sequences}
For $k\in\NN$, $k≥1$, we define a \alert{$N_k$-sequence} to be a sequence
\[\mv_j \in \Mv\quad 1 ≤ j ≤ k
\]
of $k$ independent vectors in $\Mv$,
and a sequence
\[\vv_j \in \Vv \quad 1 ≤ j ≤ k
\]
of $k$ independent vectors in $\Vv$
such that $\A\mv_j=\vv_j$ for $1≤j≤k$, $\E\mv_j=\vv_{j-1}$ for $2≤j≤k$, $\E\mv_1=0$ and $\vv_k\not\in\im\E$, which is summarized in the following diagram.
\begin{center}
\begin{tikzpicture}
\matrix(m) [exseq]
{0 & \mv_1 & \vv_1 & \cdots & \mv_k & \vv_k \not\in\im\E \\};
\EAedge{1}{2}{3}
\EAedge{4}{5}{6}
\Eedge{3}{4}
\end{tikzpicture}
\end{center}

Similarly, for $k\in\NN$, $k≥1$, we define a \alert{$L_k$-sequence} to be a sequence
\[\mv_j\in\Mv\quad 1 ≤ j ≤ k-1
\]
of $k-1$ independent vectors
 and a sequence
\[\vv_j\in\Vv\quad 1 ≤ j ≤ k
\]
of $k$ independent vectors
which fulfill the conditions summarized in the following diagram.
\begin{center}
\begin{tikzpicture}
\matrix(m) [exseq]
{\im\A \not\ni \vv_1 & \mv_1 & \cdots & \mv_{k-1} & \vv_{k} \not\in\im\E\\};
\EAedge{1}{2}{3}
\EAedge{3}{4}{5}
\end{tikzpicture}
\end{center}
\end{definition}

\todo{props of such sequences, and matrices...}

\begin{theorem}
\label{thm_NL_basis}
\autoref{thm_linear_decomposition} produces bases such that there are
$\alpha_k$ $N_k$-sequences, and $\beta_k^+$ $L_k$-sequences.
Moreover, the end vectors $\vv$ constitute a basis of $\Wv'$.
\end{theorem}

\begin{pr}
\todo{refine}
\begin{enumerate}
\item
We proceed by induction on the index.
Assume that the result holds for the reduced system $(\E',\A')$.
\item
The basis of $\Cv'$ is precisely such that $\E\mv_{k+1} = \vv_{k}$.
Besides, the basis of $\A\Cv'$ is chosen such that $\vv_{n+1} = \A\mv_{k+1}$, so each $L_k$ and $N_k$ sequence is extended with two elements, meaning that they build now $L_{k+1}$ and $N_{k+1}$ sequences.
\item
Each element $\mv$ of the basis of $\Kv'$ produces a new $N_1$ sequence $(\mv, \A\vv)$, since $\E\mv = 0$ and $\A\mv = \vv \not\in\im\E$:
\begin{center}
\begin{tikzpicture}
\matrix(m) [commdiag]
{0 & \mv & \vv \not\in\im\E\\};
\EAedge{1}{2}{3}
\end{tikzpicture}
\end{center}

The dimension of $\Kv'$ being $\alpha_1$, we produce $\alpha_1$ such sequences.
\item
Each element $\vv$ of the basis $\Zv'$ qualifies as a $L_0$ sequence, since $\vv\not\in\im\E$ and $\vv\not\in\im\A$, so
\[\im\A \not\ni \vv \not\in\im\E
.
\]
Since the dimension of $\Zv'$ is $\beta_1^+$, we produce $\beta_1^+$ such sequences.
\item
We conclude using \autoref{def_red_defects}\todo{details}
\end{enumerate}
\end{pr}

\subsection{Kronecker Decomposition}

The {Kronecker canonical form} makes use of special blocks, each of which having a variant for the matrices $\E$ and $\A$.

\newcommand{\AoneEzero}{\newcommand*{\Acol}[1]{1}\newcommand*{\Ecol}[1]{0} }
\newcommand{\EoneAzero}{\newcommand*{\Acol}[1]{0}\newcommand*{\Ecol}[1]{1} }

\begin{definition}
\label{def_bidinil}
The rectangular \alert{bidiagonal blocks} {$\mat{L}_k^{\E}$}{$\E$ part of the Kronecker $\mat{L}$-blocks} and {$\mat{L}_k^{\A}$}{$\A$ part of the Kronecker $\mat{L}$-blocks} defined by

\newcommand{\KronLBlock}{\left.\begin{bmatrix}
\Ecol{1} &  &   &  \\
\Acol{1}  & \Ecol{1} &  & \\
& \ddots & \ddots & \\
&& \Acol{1} & \Ecol{1} \\
&& & \Acol{1}
\end{bmatrix}\quad\right\} k}

\[\mat{L}_{k}^{\E} := {\EoneAzero\KronLBlock}
,
\qquad \mat{L}_k^{\A} := {\AoneEzero\KronLBlock}
.
\]

The \alert{nilpotent blocks} {$\mat{N}_k^{\E}$}{$\E$ part of the nilpotent blocks} and {$\mat{N}_k^{\A}$}{$\A$ part of the nilpotent blocks} defined by
\end{definition}

\newcommand{\Nilpotentblock}{\left. \begin{bmatrix}
	\Acol{1} & \Ecol{1} &  & & \\
	  & \Acol{1} & \Ecol{1} &  & \\
	  &   & \ddots & \ddots & \\
	& & & \Acol{1} & \Ecol{1}\\
		& & & & \Acol{1}
\end{bmatrix}\quad\right\} k}

\newcommand{\diagEA}[1]{\operatorname{diag}(
\mat{N}_{k_1}^{#1},\ldots,\mat{N}_{k_m}^{#1},
\mat{L}_{k_1}^{#1},\ldots,\mat{L}_{k_p}^{#1},
(\mat{L}_{k_1}^{#1})^{\transpose},\ldots,(\mat{L}_{k_q}^{#1})^{\transpose})}

\[ \mat{N}_k^{\E} := {\EoneAzero \Nilpotentblock}
,
\qquad \mat{N}_k^{\A} := {\AoneEzero \Nilpotentblock}
.
\]

\begin{definition}
\label{def_kronecker}
A \alert{Kronecker decomposition} of the system $(\E,\A)$ is a choice of basis of $\Mv$ and $\Vv$ such that $\E$ and $\A$ are decomposed in blocks of the same size
\[\E = \begin{bmatrix}
	 \Id & 0 \\ 0 & \overline{\E}
\end{bmatrix}, \qquad \A = \begin{bmatrix}
	\mat{J} & 0 \\ 0 & \overline{\A}
\end{bmatrix}
,
\]
where $\mat{J}$ is a diagonal block of Jordan blocks, and $\overline{E}$ and $\overline{\A}$ are in diagonal block form
\[\begin{aligned}
\overline{\E} &= \diagEA{\E}
, \\
\overline{\A} &= \diagEA{\A}
,
\end{aligned} \]
where the blocks of $\E$ and $\A$ have the same size.
\end{definition}

\begin{theorem}
\label{thm_Kronecker}
A decomposition with defects $α$, $β^+$, $β^-$, produces a Kronecker decomposition which for all integer $k≥1$ contains
\begin{itemize}
	\item $α_k$ block of type $\mat{N}_{k}$,
	\item $β_k^+$ blocks of type $\mat{L}_k$,
	\item $β_k^-$ blocks of type $\mat{L}_k^{\transpose}$.
\end{itemize}
\end{theorem}

\begin{pr}
Recall the definition of $L_k$ and $N_k$ sequences in \autoref{def_NL_sequences}.
By regrouping the elements of an $L_k$ sequence, one obtains a representation of $\E$ and $\A$ as a $\mat{L}_k$-block, and similarly, by regrouping the elements of a $N_k$-sequence, one obtains a $\mat{N}_k$ block.
Applying \autoref{thm_NL_basis}, and regrouping the basis elements stemming from the sequences $L_k$ and $N_k$ we obtain $α_k$ $\mat{N}_k$-blocks and $β^+_k$ $\mat{L}_k$-blocks, for $k≥1$.

Now the basis on the sub-block $\Mv\red{∞}$, $\Vv\red{∞}$ are obtained by transposing the decomposition given by \autoref{thm_linear_decomposition}. 
Using the previous step and \autoref{prop_alpha_reds_zero}, we obtain $β^-_k$ transpose of $\mat{L}_k$-blocks, for $k≥1$.
\end{pr}

\begin{remark}
It is remarkable that the decomposition obtained in \autoref{thm_full_linear_decomposition} produces basis vectors which are \emph{the same} as for a Kronecker decomposition, only ordered differently.
The necessary permutations may be visualised on \autoref{fig_kronecker}.
\end{remark}

\matrixfigure{fig_kronecker}{kronecker_fig}{
An illustration of the defects $α$, $β^+$ and $β^-$ and of the Kronecker decomposition described in \autoref{thm_Kronecker}.
The difference of size of the squares is exactly given by the defects $α$, $β^+$ and $β^-$.
The dark squares bearing the number $j$ represent all the nilpotent blocks $\mat{N}_j$; there are $α_j$ such blocks.
The light squares in the lower-right part bearing the number $j$ represent the $\mat{L}$-blocks $\mat{L}_j$.
There are $β_j^+$ such blocks.
The light squares in the upper-left part bearing the number $j$ represent the $\mat{L}$-blocks $\mat{L}_j^{\transpose}$.
There are $β_j^-$ such blocks.
This figure also allows to check the formulae of \autoref{prop_rel_defect_dims}.
}




\subsection{Conjugate Decomposition}

We may now show the relation between the defects $α$, $β^+$ and $β^-$ of a system $(\E,\A)$ and the defects of the adjoint system $(\E^*, \A^*)$.
It turns out that the constraint defects are the same and that the observation defects $β^+$ and the control defects $β^-$ are just switched.
This fact would have been very difficult to prove from the results of \autoref{sec_defects} alone, so we need the full power of \autoref{thm_full_linear_decomposition} and of its consequence, \autoref{thm_Kronecker}.

\begin{theorem}
\label{thm_conjugate_defects}
	The {conjugate decomposition} switches the defects $β^+$ and $β^-$, i.e., it produces the  defects
	\[\begin{aligned}
		α(\E^*,\A^*) &= α(\E,\A)
		,\\
		{β^+}(\E^*,\A^*) &= β^-(\E,\A)
		,\\
		{β^-}(\E^*,\A^*) &= β^+(\E,\A)
		.
	\end{aligned}\]
\end{theorem}
\begin{pr}
It is a consequence of \autoref{thm_Kronecker}, for when putting the system in Kronecker form \autoref{def_kronecker} and transposing, the system is still in Kronecker form, but the bidiagonal blocks $\mat{L}$ and $(\mat{L})^{\transpose}$ are switched.
\end{pr}

\subsection{Weierstraß decomposition}

In the case of regular pencils (see \autoref{sec_regular_pencil}), the Kronecker decomposition is called the \emph{Weierstraß decomposition} (\cite{Weierstrass,Gantmacher}) and is such that $\E$ and $\A$ take the matrix representation
\[\E = \begin{bmatrix}
	\Id & 0 \\ 0 & \mat{N}
\end{bmatrix}
,
\qquad \A = \begin{bmatrix}
	\mat{C} & 0 \\ 0 & \Id
\end{bmatrix}
,
\]
where $\mat{C}$ may be in Jordan normal form and $\mat{N}$ is a block diagonal matrix of blocks of type $\mat{N}_k^E$. 

The matrix block {$\mat{N}$}{Matrix of nilpotent blocks} is a diagonal block matrix
\[\mat{N} = \mathrm{diag}\bigl(\mat{N}^E_{k_1}(0),\mat{N}^E_{k_2}(0),\ldots,\mat{N}^E_{k_m}(0)\bigr)
\]
where the blocks $(\mat{N}^E_{k_1}(0)$ are the nilpotent blocks defined in \autoref{def_bidinil}.

\begin{corollary}
The Weierstraß decomposition is such that for any integer $k≥1$ it contains
$α_k$ blocks $\mat{N}_k^{\E}$.
\end{corollary}
\begin{pr}
It is just a special case of \autoref{thm_Kronecker} using \autoref{prop_reg_pencil}.
\end{pr}

\section{Conclusions}

We have defined the notion of \emph{defects} and have related them to existing concepts, such as the regular pencil condition, the dimension of the reduced subspaces, or the notion of strangeness.
We also showed how the defects define a normal form, and how that normal form relates to the existing one of Kronecker.

Note that some results, as \autoref{thm_conjugate_defects}, would be difficult to prove without using the canonical form.
Nevertheless, we tried to wring the most out of the invariant objects defined in \autoref{sec_defects}.

The advantage of such an approach is that it is extensible to nearby cases such as the parameter dependent case, or the infinite dimensional case (see \cite{linpdae}).

\newcommand{\includebibliography}[1]{
\bibliographystyle{abbrv}
\bibliography{../#1}
}

\includebibliography{IDE}

\end{document}